%% file: main.tex
\documentclass [10pt,reqno]{amsart}
\usepackage {hannah}
\usepackage {color,graphicx,tikz}
\usepackage {multirow}
\usepackage [utf8]{inputenc}
\usepackage [bookmarksnumbered=true,colorlinks=true]{hyperref}

\newcommand {\calC}{{\mathcal C}}
\newcommand {\calM}{{\mathcal M}}
\newcommand {\calN}{{\mathcal N}}

\newcommand {\PP}{{\mathbb P}}

\newcommand {\RR}{{\mathbb R}}
\newcommand {\ZZ}{{\mathbb Z}}

\newcommand {\rdim}{\textnormal{rdim}}
\newcommand {\RH}{\textnormal{RH}}

\newcommand{\mm}[2]{\calM_{0}(#1,#2)}

\DeclareMathOperator {\ev}{ev}
\DeclareMathOperator {\ft}{f\/t}
\DeclareMathOperator {\val}{val}

\newcommand {\dunion}{\,\mbox {\raisebox{0.25ex}{$\cdot$} \kern-1.83ex $\cup$}
  \,}

\title [Tropical moduli spaces of stable maps to a curve]{Tropical moduli
  spaces of stable maps to a curve}
\author {Andreas Gathmann, Hannah Markwig, and Dennis Ochse}

\address {Andreas Gathmann, Fachbereich Mathematik, Technische Universit\"at
  Kaiserslautern, Postfach 3049, 67653 Kaiserslautern, Germany}
\email {andreas@mathematik.uni-kl.de}

\address {Hannah Markwig, Eberhard Karls Universit\"at T\"ubingen, Fachbereich Mathematik, Institut f\"ur Geometrie, Auf der Morgenstelle 10, 72076 T\"ubingen, Germany}
\email{hannah@math.uni-tuebingen.de}

\address{Dennis Ochse, Fachbereich Mathematik, Technische Universit\"at
  Kaiserslautern, Postfach 3049, 67653 Kaiserslautern, Germany}
\email{ochse@mathematik.uni-kl.de}


\newcommand {\refx}[1]{\ref*{#1}}

\begin {document}

\begin {abstract}
  We construct moduli spaces of rational covers of an arbitrary smooth tropical
  curve in $\R^r$ as tropical varieties. They are contained in the balanced fan
  parametrizing tropical stable maps of the appropriate degree to $\R^r$. The
  weights of the top-dimensional polyhedra are given in terms of certain
  lattice indices and local Hurwitz numbers. 
\end {abstract}

\maketitle


\input{introduction.tex}
\input{preliminaries.tex}

\input{gluing.tex}

\input{curves.tex}

\bibliographystyle{amsalpha}
\bibliography{bibliographie}

\end {document}

%% file: introduction.tex
\section{Introduction} \label{intro}

Tropical enumerative geometry has developed from interesting applications
following so-called correspondence theorems which settle the equality of
certain enumerative numbers in algebraic geometry to their tropical
counterparts \cite{Mi03}. There is an ongoing effort to put the striking
similarities between algebro-geometric and tropical enumerative geometry onto a
more solid ground.

Modern enumerative algebraic geometry is based on the moduli spaces $\overline
M_{g,n}(X,\beta)$ of $n$-pointed stable maps of genus $g$ and class $ \beta $
to a smooth projective variety $X$ \cite{FP97}, together with their virtual
fundamental classes \cite{BF97,Behrend} that resolve the issues arising when
these spaces are not of the expected dimension. Hence a key ingredient for the
further development of tropical enumerative geometry is the construction of
tropical analogues of these concepts. If $ g=0 $ and $X$ is a toric variety,
corresponding to rational tropical curves in $\R^r$, such tropical spaces have
been constructed as balanced fans in \cite{GKM07}. In this case, ideas relating
to virtual fundamental classes are not needed, and the intersection theory of
the resulting spaces recovers the correspondence theorems for rational tropical
curves in $\R^r$ \cite{Gro14}.

For more general target spaces, we run into the same problems as in algebraic
geometry: the naively defined spaces of tropical curves in a tropical variety
are usually not of the expected dimension, maybe not even pure-dimensional.
However, as there is no general theory of virtual fundamental classes in
tropical geometry yet, the tropical approach to this problem is different:
right from the start we have to construct the moduli spaces as balanced
polyhedral complexes of the expected dimension --- which necessarily means that
they are \emph{not} just the spaces of maps from a tropical curve to the given
target. From an algebro-geometric point of view, one could say that this
constructs the moduli space and its virtual fundamental class at the same time,
with the additional benefit that (in accordance with the general philosophy of
tropical intersection theory) we actually obtain a virtual \emph{cycle} and not
just a \emph{cycle class}.

A general approach how this idea might be realized has been presented in
\cite{GO14}. Here, we will restrict ourselves to the case when $ g=0 $ and the
target is a smooth (rational) tropical curve $L$ in $ \R^r $. The resulting
moduli spaces $ \calM_{0,n}(L,\Sigma) $ (where $ \Sigma $ is a degree of
tropical curves as in definition \ref{def-curve}) then describe rational covers
of a rational smooth tropical curve.

Tropical covers and tropical Hurwitz numbers (i.e.\ enumerative numbers
counting covers with prescribed properties \cite{BBM10, CJM10}) are useful
e.g.\ for the study of the structural behavior of Hurwitz numbers
\cite{cjm:wcfdhn} and in the tropical enumeration of Zeuthen numbers
\cite{BBM11}. Spaces of tropical (admissible) covers have been studied in
\cite{CMR14} as tropicalizations of corresponding algebro-geometric spaces, in
terms of a tropicalization map on the Berkovich analytification. The space of tropical covers of $\R$ has been described in \cite{CMR14b} as tropicalization of the open part of a suitable space of relative stable maps (whose compactification is then realized as a tropical compactification defined by the tropical moduli space). The present
work complements this point of view by fixing a rational smooth tropical curve
$L\subset \R^r$, restricting to genus $0$ covers, and embedding the abstract
polyhedral subcomplex of the abstract cone complex described in \cite{CMR14} as
a balanced polyhedral subcomplex. In this way, we make these moduli spaces
accessible to the current state of the art of tropical intersection theory.

As mentioned above, to construct the moduli spaces $\mathcal M_{0,n}(L,\Sigma)$
we cannot just take the subset of $\mathcal M_{0,n}(\R^r,\Sigma)$ consisting of
all stable maps whose image lies in $L$, as this would yield a non-pure
subcomplex with strata of too big dimension. Instead, we have to incorporate
the so-called Riemann-Hurwitz condition (see definition
\ref{def-RH}), which implies the algebraic realizability of the corresponding
maps. For an example, let $ L \subset \R^2 $ be the standard tropical line, let
$ \Sigma $ be the degree consisting of the directions $
(-1,0),(-1,0),(0,-1),(0,-1),(2,2) $, and set $ n=0 $. A fan curve in $L$ of
this degree --- in fact representing the origin of the fan $ \mathcal
M_{0,0}(L,\Sigma) $ --- is shown in picture (a) below. It is given by a map
from an abstract star curve with $5$ ends to $L$, with the directions and
weights on the ends as indicated in the picture.

\begin {center} \input {pics/intro} \end {center}

Possible resolutions of this curve in $L$ are shown in (b), (c), and (d).
However, case (b) is excluded in $ \calM_{0,0}(L,\Sigma) $ as its central
vertex violates the Riemann-Hurwitz condition: it would correspond to an
algebraic degree-$2$ cover of the projective line by itself with three
ramification points of order $2$, which does not exist. In contrast, the
combinatorial types (c) and (d) are allowed, and represent two rays in $
\calM_{0,0}(L,\Sigma) $ since they describe $1$-dimensional families of curves.
They both have a similar type obtained by symmetry: in (c) the bounded
weight-$2$ edge could also be on the horizontal edge of $L$, and in (d) there
are two choices how to group the weight-$1$ ends. In total, this means that $
\calM_{0,0}(L,\Sigma) $ is a $1$-dimensional fan with four rays. The weights
that we will construct on these rays incorporate the triple Hurwitz numbers
corresponding to the local degrees of the maps at each point mapping to the
vertex of $L$; they all turn out to be $1$ here. In this example, it is then
easy to check explicitly that $ \calM_{0,0}(L,\Sigma) \subset
\calM_{0,0}(\R^2,\Sigma) \cong \calM_{0,5} \times \R^2 $ is indeed balanced.
Our main result on the moduli spaces $\mathcal M_{0,n}(L,\Sigma)$ is that this
construction works in general:

\begin{theorem} \label{thm1}
  Let $L$ be a smooth tropical curve in $\R^r$ and $\Sigma$ a degree of
  tropical stable maps to $L$ (see definitions \ref{def-curve} and \ref{def-ml}). Then the space
  $\mathcal M_{0,n}(L,\Sigma)$ (with weights defined in terms of local
  Hurwitz numbers) is a balanced weighted polyhedral subcomplex of $
  \mathcal M_{0,n}(\R^r,\Sigma)$ of pure dimension
    \[ |\Sigma|- \deg(\Sigma)\cdot \Big(\sum_{W \in L} (\val(W)-2) \Big)-2. \]
\end{theorem}

We expect that $\mathcal M_{0,n}(L,\Sigma)$ is in fact the tropicalization of
(relevant parts) of the corresponding algebro-geometric moduli space.

Theorem \ref{thm1} is proved in two major steps: the first being the treatment
of $1$-dimensional moduli spaces of the form above (see theorem
\ref{thm-local-balancing}), and the second the generalization to arbitrary
dimension. For the generalization to arbitrary dimension, we use a general
gluing construction for tropical moduli spaces which was developed by the first
and last author in \cite{GO14} and has further applications to other target
spaces.
 
This paper is organized as follows. In section \ref{sec-prelim} we review the
necessary preliminaries. The tropical moduli spaces $\mathcal
M_{0,n}(L,\Sigma)$ are then defined in section \ref{sec-m0nldelta}. More
precisely, we define their structure as a polyhedral subcomplex of $\mathcal
M_{0,n}(\R^r,\Sigma)$ in subsection \ref{subsec-ascomplex}, and the weights of
their maximal cells in subsection \ref{subsec-gluing}. The definition of the
weights relies on the gluing construction of \cite{GO14}, which we therefore
review in subsection \ref{subsec-gluing}, together with the main result of
\cite{GO14} allowing a gluing construction of tropical moduli spaces under some
requirements. In our case, these requirements are satisfied if all
one-dimensional tropical moduli spaces $\mathcal M_{0,n}(L,\Sigma)$ are
balanced fans. We prove this fact in section \ref{sec-1dim} (see theorem
\ref{thm-local-balancing}). Theorem \ref{thm1} is then an immediate consequence
of the foundational work on the gluing construction in \cite{GO14}. 

\subsection{Acknowledgments}

We would like to thank Erwan Brugall\'{e}, Renzo Cavalieri, Simon Hampe, and Diane Maclagan
for helpful discussions.

This work would not have been possible without extensive computations of example classes which enabled us to establish and prove conjectures about polyhedra and their weights in our moduli spaces. We used the polymake-extension a-tint \cite{polymake, Ham12} and GAP \cite{gap4}.

Part of this work was accomplished at the Mittag-Leffler Institute in
Stockholm, during the semester program in spring 2011 on Algebraic Geometry
with a View towards Applications. The authors would like to thank the institute
for hospitality.

The first and second author were partially funded by DFG grant GA 636/4-2 resp.\ MA 4797/3-2, as part of the DFG Priority Program 1489.

We thank an anonymous referee for helpful comments on an earlier version of this paper.

%% file: pics/intro.tex
\begin{picture}(0,0)%
\includegraphics{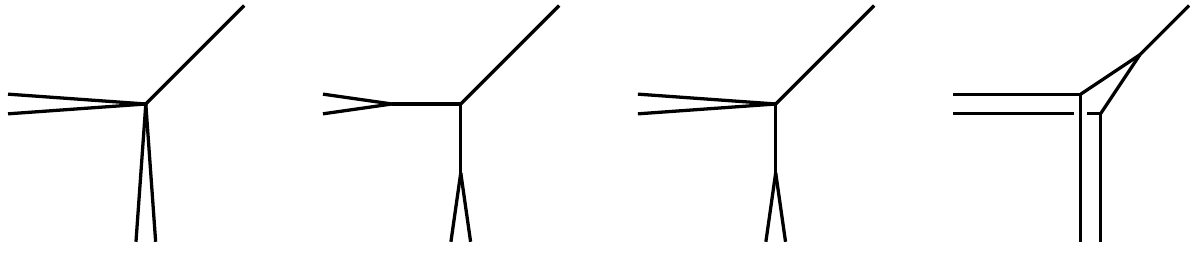}%
\end{picture}%
\setlength{\unitlength}{4144sp}%
\begingroup\makeatletter\ifx\SetFigFont\undefined%
\gdef\SetFigFont#1#2#3#4#5{%
  \reset@font\fontsize{#1}{#2pt}%
  \fontfamily{#3}\fontseries{#4}\fontshape{#5}%
  \selectfont}%
\fi\endgroup%
\begin{picture}(5459,1156)(1539,-430)
\put(2611,524){\makebox(0,0)[lb]{\smash{{\SetFigFont{8}{9.6}{\familydefault}{\mddefault}{\updefault}{\color[rgb]{0,0,0}$2$}%
}}}}
\put(1666,344){\makebox(0,0)[b]{\smash{{\SetFigFont{8}{9.6}{\familydefault}{\mddefault}{\updefault}{\color[rgb]{0,0,0}$1$}%
}}}}
\put(1666, 74){\makebox(0,0)[b]{\smash{{\SetFigFont{8}{9.6}{\familydefault}{\mddefault}{\updefault}{\color[rgb]{0,0,0}$1$}%
}}}}
\put(2116,-331){\makebox(0,0)[rb]{\smash{{\SetFigFont{8}{9.6}{\familydefault}{\mddefault}{\updefault}{\color[rgb]{0,0,0}$1$}%
}}}}
\put(2296,-331){\makebox(0,0)[lb]{\smash{{\SetFigFont{8}{9.6}{\familydefault}{\mddefault}{\updefault}{\color[rgb]{0,0,0}$1$}%
}}}}
\put(4051,524){\makebox(0,0)[lb]{\smash{{\SetFigFont{8}{9.6}{\familydefault}{\mddefault}{\updefault}{\color[rgb]{0,0,0}$2$}%
}}}}
\put(3106,344){\makebox(0,0)[b]{\smash{{\SetFigFont{8}{9.6}{\familydefault}{\mddefault}{\updefault}{\color[rgb]{0,0,0}$1$}%
}}}}
\put(3106, 74){\makebox(0,0)[b]{\smash{{\SetFigFont{8}{9.6}{\familydefault}{\mddefault}{\updefault}{\color[rgb]{0,0,0}$1$}%
}}}}
\put(3489,299){\makebox(0,0)[b]{\smash{{\SetFigFont{8}{9.6}{\familydefault}{\mddefault}{\updefault}{\color[rgb]{0,0,0}$2$}%
}}}}
\put(3678, 51){\makebox(0,0)[lb]{\smash{{\SetFigFont{8}{9.6}{\familydefault}{\mddefault}{\updefault}{\color[rgb]{0,0,0}$2$}%
}}}}
\put(3570,-331){\makebox(0,0)[rb]{\smash{{\SetFigFont{8}{9.6}{\familydefault}{\mddefault}{\updefault}{\color[rgb]{0,0,0}$1$}%
}}}}
\put(3723,-331){\makebox(0,0)[lb]{\smash{{\SetFigFont{8}{9.6}{\familydefault}{\mddefault}{\updefault}{\color[rgb]{0,0,0}$1$}%
}}}}
\put(5491,524){\makebox(0,0)[lb]{\smash{{\SetFigFont{8}{9.6}{\familydefault}{\mddefault}{\updefault}{\color[rgb]{0,0,0}$2$}%
}}}}
\put(4546,344){\makebox(0,0)[b]{\smash{{\SetFigFont{8}{9.6}{\familydefault}{\mddefault}{\updefault}{\color[rgb]{0,0,0}$1$}%
}}}}
\put(4546, 74){\makebox(0,0)[b]{\smash{{\SetFigFont{8}{9.6}{\familydefault}{\mddefault}{\updefault}{\color[rgb]{0,0,0}$1$}%
}}}}
\put(5118, 51){\makebox(0,0)[lb]{\smash{{\SetFigFont{8}{9.6}{\familydefault}{\mddefault}{\updefault}{\color[rgb]{0,0,0}$2$}%
}}}}
\put(5010,-331){\makebox(0,0)[rb]{\smash{{\SetFigFont{8}{9.6}{\familydefault}{\mddefault}{\updefault}{\color[rgb]{0,0,0}$1$}%
}}}}
\put(5163,-331){\makebox(0,0)[lb]{\smash{{\SetFigFont{8}{9.6}{\familydefault}{\mddefault}{\updefault}{\color[rgb]{0,0,0}$1$}%
}}}}
\put(6931,524){\makebox(0,0)[lb]{\smash{{\SetFigFont{8}{9.6}{\familydefault}{\mddefault}{\updefault}{\color[rgb]{0,0,0}$2$}%
}}}}
\put(5986,344){\makebox(0,0)[b]{\smash{{\SetFigFont{8}{9.6}{\familydefault}{\mddefault}{\updefault}{\color[rgb]{0,0,0}$1$}%
}}}}
\put(5986, 74){\makebox(0,0)[b]{\smash{{\SetFigFont{8}{9.6}{\familydefault}{\mddefault}{\updefault}{\color[rgb]{0,0,0}$1$}%
}}}}
\put(6450,-331){\makebox(0,0)[rb]{\smash{{\SetFigFont{8}{9.6}{\familydefault}{\mddefault}{\updefault}{\color[rgb]{0,0,0}$1$}%
}}}}
\put(6603,-331){\makebox(0,0)[lb]{\smash{{\SetFigFont{8}{9.6}{\familydefault}{\mddefault}{\updefault}{\color[rgb]{0,0,0}$1$}%
}}}}
\put(5874,-354){\makebox(0,0)[lb]{\smash{{\SetFigFont{10}{12.0}{\familydefault}{\mddefault}{\updefault}{\color[rgb]{0,0,0}(d)}%
}}}}
\put(1554,-354){\makebox(0,0)[lb]{\smash{{\SetFigFont{10}{12.0}{\familydefault}{\mddefault}{\updefault}{\color[rgb]{0,0,0}(a)}%
}}}}
\put(2994,-354){\makebox(0,0)[lb]{\smash{{\SetFigFont{10}{12.0}{\familydefault}{\mddefault}{\updefault}{\color[rgb]{0,0,0}(b)}%
}}}}
\put(4434,-354){\makebox(0,0)[lb]{\smash{{\SetFigFont{10}{12.0}{\familydefault}{\mddefault}{\updefault}{\color[rgb]{0,0,0}(c)}%
}}}}
\put(6661,209){\makebox(0,0)[lb]{\smash{{\SetFigFont{8}{9.6}{\familydefault}{\mddefault}{\updefault}{\color[rgb]{0,0,0}$1$}%
}}}}
\put(6571,389){\makebox(0,0)[rb]{\smash{{\SetFigFont{8}{9.6}{\familydefault}{\mddefault}{\updefault}{\color[rgb]{0,0,0}$1$}%
}}}}
\end{picture}%

%% file: preliminaries.tex
\section{Preliminaries} \label{sec-prelim}

\subsection{Background on tropical varieties and intersection theory}
  \label{subsec-X}

To fix notation, we quickly review notions of tropical intersection theory.
Some of our constructions involve partially open versions of tropical
varieties, i.e.\ varieties containing polyhedra that are open at some faces. We
adapt the usual conventions to this situation. For a more detailed survey of
the relevant preliminaries, see \cite[section \ref{modulx:subsec-X}]{GO14}.

We let $\Lambda$ be a lattice in an $r$-dimensional real vector space $V$. A
\emph{(partially open) (rational) polyhedron} in $V$ is a finite intersection
of (open or) closed affine half-spaces given by (strict or) non-strict
inequalities whose linear parts are given by elements in the dual of $\Lambda$.
We denote by $V_\sigma$ the linear space obtained by shifting the affine span
of $\sigma$ to the origin and define $ \Lambda_\sigma := V_\sigma \cap
\Lambda$. A \emph{face} $ \tau \leq \sigma $ (also written as $\tau<\sigma$ if
it is proper) is a non-empty subset of $ \sigma $ that can be obtained by
changing some of the defining non-strict inequalities into equalities. If $
\dim \tau = \dim \sigma -1 $ we call $ \tau $ a \emph{facet} of $ \sigma $. In
this case we denote by $ u_{\sigma/\tau} \in \Lambda_\sigma / \Lambda_\tau $
the \emph{primitive normal vector} of $\sigma$ relative to $\tau$, i.e.\ the
unique generator of $ \Lambda_\sigma / \Lambda_\tau $ lying in the half-line of
$\sigma$ in $ V_\sigma / V_\tau \cong \RR $. The well-known notion of a
(pure-dimensional) \emph{weighted polyhedral complex} $X$ (formed by cells
$\sigma$ as above, and with integer weights on maximal cells), its dimension
and support are easily adapted to the case of partially open polyhedral
complexes. Such a (partially open) weighted polyhedral complex $ (X,\omega) $
is called a \emph{(partially open) tropical variety} (or \emph{cycle}, if
negative weights occur) if it satisfies the balancing condition, i.e.\ for each
cell $\tau$ of codimension $1$ we have
\begin{displaymath}
  \sum_{\sigma:\sigma>\tau} \omega(\sigma) \cdot u_{\sigma/\tau} = 0
    \quad \in V / V_\tau.
\end{displaymath}
For intersection-theoretic purposes, the exact polyhedral complex structure is
often not important, and we fix it only up to refinements respecting the weights.

\begin{example}[Smooth curves] \label{ex-curves}
  Let $V=\RR^q$. We let $L_1^q$ denote the $1$-dimensional tropical variety
  containing the origin and rays spanned by $-e_i$ (where $e_i$ denotes the
  canonical basis vectors) and $-e_0:=\sum e_i$, with all weights one. This is
  the tropicalization of a general line over the Puiseux series with constant
  coefficient equations \cite[proposition 2.5 and theorem 4.1]{FS05}. A
  one-dimensional tropical variety $L\subset \R^r$ with all weights one is
  called a \emph{rational smooth curve} if its underlying polyhedral complex
  is rational (i.e.\ combinatorially a tree), and if it locally at each vertex
  equals $L_1^q$ up to a \emph{unimodular transformation}, i.e.\ up to an
  isomorphism of vector spaces which is also an isomorphism of the underlying
  lattices \cite{All09}.
\end{example}

Some of our constructions involve \emph{quotients} $X/W$ of partially open
tropical varieties $X$ by a \emph{lineality space} $W$. We say that a vector
subspace $W$ of $V$ is a lineality space for $X$ if for all $ \sigma \in X $
and $ x \in \sigma $ the intersection $ \sigma \cap (x+L) $ is open in $ x+L $
and equal to $ |X| \cap (x+L) $. Note that for the case of a closed polyhedral
complex this generalizes the usual notion of a lineality space (which is
commonly the maximal subspace with this property). For more details on such
quotients, see \cite[section \refx{modulx:subsec:quotients}]{GO14}.

A \emph{morphism} between (partially open) tropical cycles $X$ and $Y$ is a map
$ f: |X| \to |Y| $ which is locally affine linear, with the linear part induced
by a map between the underlying lattices \cite[definition 7.1]{AR07}. A
\emph{rational function} on a tropical variety $X$ is a continuous function
$\varphi : |X| \to \RR$ that is affine linear on each cell, and whose linear
part is integer, i.e.\ in the dual of the lattice. We associate a
\emph{divisor} $\varphi \cdot X$ to a rational function; a cycle of codimension
$1$ in $X$ support on the cells at which $ \varphi $ is not locally linear
\cite[construction 3.3]{AR07}. Multiple intersection products $\varphi_1 \cdot
\; \cdots \; \cdot \varphi_m \cdot X$ are commutative by \cite[proposition
3.7]{AR07}.

\begin{remark}[Weights of intersections as lattice indices]
    \label{rem-latticeindex}
  Often, the weight of a cell of a multiple intersection product can be
  computed locally in terms of a lattice index. To do this, we write locally $
  \varphi_i = \max\{h_i,0\}$ for linearly independent integer linear functions
  $ h_1, \ldots, h_m $, and let $H$ be a matrix representing the integer linear
  map $\Lambda \rightarrow \ZZ^m: x \mapsto (h_1(x), \ldots, h_m(x))$. Then the
  local weight of $\varphi_1 \cdot \ldots \cdot \varphi_m \cdot X$ equals the
  greatest common divisor of the maximal minors of $H$
  \cite[lemma 5.1]{MR08}.
\end{remark}

Rational functions can be pulled back along a morphism $f:X\rightarrow Y$ to
rational functions $f^*(\varphi) = \varphi \circ f$ on $X$. We can push forward
a subvariety $Z$ of $X$ to a subvariety $f_*(Z)$ of $Y\subset \Lambda'
\otimes_\ZZ \RR$ \cite[proposition 4.6 and corollary 7.4]{AR07}: For suitable
refinements of the polyhedral structures of $X$ and $Y$, we obtain $f(\sigma)
\in Y$ for all $\sigma \in X$, and define the weight of the push-forward to be 
  $$ \omega_{f_*(Z)} (\sigma') :=
       \sum_{ \sigma}
       \omega_X (\sigma) \cdot |\Lambda'_{\sigma'}/f(\Lambda_\sigma)|, $$
where the sum goes over all top-dimensional cells $\sigma\in Z$ with $f(\sigma)
= \sigma' $. In the partially open case, we will restrict ourselves to
injective morphisms in order to avoid problems with overlapping cells with
different boundary behavior.

\subsection{Tropical moduli spaces of curves} \label{subsec-m0nrd}

An \emph{(abstract) $N$-marked rational tropical curve} is a tuple $
(\Gamma,x_1,\dots,x_N) $, where $\Gamma$ is a metric tree with $N$ unbounded
edges labeled $ x_1,\dots,x_N $ (also called \emph{marked ends}) that have
infinite length, and such that the valence of each vertex is at least $3$. The
set of all $N$-marked tropical curves is denoted $ \calM_{0,N}$. It follows
from \cite[theorem 3.4]{SS04a}, \cite[section 2]{Mi07}, or \cite[theorem
3.7]{GKM07} that $\calM_{0,N}$ can be embedded as a tropical variety via the
distance map, more precisely, as a balanced, simplicial, $(N-3)$-dimensional
fan whose top-dimensional cones all have weight one. The distance map sends a
tropical curve to the vector of distances of its ends in $\R^{\binom{N}{2}}$.
We mod out an $N$-dimensional lineality space $U_N$, identifying vectors
corresponding to trees whose metrics only differ on the ends. For a tree with
only one bounded edge of length one, the ends with markings $I\subset
\{1,\ldots,N\}$, $1<|I|<N-1$, on one side and the ends with markings $I^c$ on
the other, we denote the equivalence class of its image under the distance map
in $\R^{\binom{N}{2}}/U_N$ by $v_I$. The vectors $v_I$ generate the rays of
$\calM_{0,N}$ and the lattice we fix for $\R^{\binom{N}{2}}/U_N$.

For local computations, we sometimes use a finite index set $I$ instead of $
\{1,\dots,N\} $ as labels for the markings, and denote the corresponding moduli
spaces by $\calM_{0,I}$. Also, we can modify the definition above by assigning
bounded lengths in $ \R_{>0} $ to the ends, corresponding to not taking the
quotient by $ U_N $. In this case we obtain a partially open moduli space which
we will denote by $ \calM'_{0,N} $. There is then a map $ \calM'_{0,N} \to
\calM_{0,N} $ forgetting the lengths of the bounded ends, which is just the
quotient by $U_N$.

For every subset $I\subset \{1,\dots,N\} $ of cardinality at least three, there
is a \emph{forgetful map} $\ft_I:\calM_{0,N}\to \calM_{0,|I|}$ which maps
$(\Gamma,x_1,\dots,x_N)$ to the tree where we remove all ends $ x_i $ with
labels $ i \notin I $ (and possibly straighten $2$-valent vertices). Forgetful
maps are morphisms by \cite[proposition 3.9]{GKM07}. In coordinates, we project
to distances of ends in $I$.

\begin{lemma} \label{lem-forget}
  A vector $x$ in $\R^{\binom{N}{2}}/U_N$ is zero if and only if $\ft_I(x)=0$
  for all $ I\subset \{1,\dots,N\} $ with $|I|=4$.
\end{lemma}

\begin{proof}
  As $\ft_I$ is linear, the ``only if'' direction is obvious. For the other
  direction, denote the standard basis vectors of $\RR^{\binom{N}{2}}$ by
  $e_{ij}$ for $ i<j $. Let $\tilde{x} = \sum_{i<j}\lambda_{ij} e_{ij} \in
  \RR^{\binom{N}{2}}$ be a representative of $x$. For any $I$ with $ |I|=4 $,
  the assumption $ \ft_I(x)=0 $ means that the projection $
  \sum_{i,j \in I; i<j} \lambda_{ij}e_{ij}$ is in $U_4$. By definition of $
  U_4 $, it follows that there is a vector $ \mu \in \RR^I $ such that $
  \lambda_{ij}=\mu_i+\mu_j $ for all $ i<j $ in $I$, and thus that $
  \lambda_{ik}+\lambda_{jl}=\lambda_{ij}+\lambda_{kl} $ if $ I= \{i,j,k,l\} $.

  But this means that for all $ i=1,\dots,N $ the assignment
    \[ \lambda_i := \frac{1}{2}(\lambda_{ij}+\lambda_{ik}-\lambda_{jk})
       \quad \text {for arbitrary $ j,k \neq i $} \]
  is well-defined, because if $m$ is another index we have
  \begin{align*}
    \frac{1}{2}(\lambda_{ij}+\lambda_{ik}-\lambda_{jk})
      &= \frac{1}{2}(\lambda_{im}+\lambda_{ik}-\lambda_{mk})
        +\frac{1}{2}(\lambda_{ij}-\lambda_{im}+\lambda_{mk}-\lambda_{jk}) \\
      &= \frac{1}{2}(\lambda_{im}+\lambda_{ik}-\lambda_{mk}).
  \end{align*}
  As the definition of $ \lambda_i $ also implies that $ \lambda_{ij} =
  \lambda_i+\lambda_j $ for all $ i<j $, we conclude that $ \tilde{x}\in U_N $,
  and hence $x=0$.
\end{proof}

\begin{definition}[Tropical stable maps] \label{def-curve}
  Let $ n \in \N $ and $ N \ge n $. Consider a tuple $ (\Gamma,x_1,\dots,x_N,h)
  $, where $ (\Gamma,x_1,\dots,x_N) $ is an $N$-marked abstract rational tropical curve
  and $ h: \Gamma \to \RR^r $ is a continuous map that is integer linear on
  each edge. For an edge $e$ starting at a vertex $V$ of $ \Gamma $, we denote
  the tangent vector of $ h|_e $ at $V$ by $ v(e,V) \in \ZZ^r $ and call it the
  \emph{direction} of $e$ at $V$. If $e$ is an end and $V$ its only neighboring
  vertex we write $ v(e,V) $ also as $v(e)$ for simplicity.

  We say that $ (\Gamma,x_1,\dots,x_N,h) $ is an \emph{$n$-marked (rational) tropical
  stable map} to $ \R^r $, also called a \emph {(parametrized) $n$-marked curve
  in} $ \RR^r $ \cite[definition 4.1]{GKM07}, if
  \begin{itemize}
  \item $h$ satisfies the balancing condition $ \sum_{e \ni V} v(e,V) = 0 $ at
    each vertex $V$ of $ \Gamma $;
  \item $ v(x_i)=0 $ for $ i=1,\dots,n $ (i.e.\ each of the first $n$ ends is
    contracted by $h$), whereas $ v(x_i) \neq 0 $ for $ i>n $ (i.e.\ the
    remaining $ N-n $ ends are ``non-contracted ends'').
  \end{itemize}
  Two $n$-marked tropical stable maps $ (\Gamma,x_1,\dots,x_N,h) $ and $
  (\tilde \Gamma, \tilde x_1,\dots,\tilde x_N,\tilde h)$ in $ \RR^r $ are
  isomorphic (and will from now on be identified) if there is an isomorphism
  $\varphi$ of the underlying $N$-marked abstract curves such that $
  \tilde h \circ \varphi = h $.

  The \emph{degree} of an $n$-marked tropical stable map is the $N$-tuple
    $$ \Sigma = (v(x_{1}),\dots,v(x_N)) \in(\ZZ^r )^{N} $$
  of directions of its ends, including the zero directions at the first $n$
  ends. Its \emph{combinatorial type} is given by the data of the combinatorial
  type of the underlying abstract marked tropical curve $ (\Gamma,x_1,\dots,
  x_N) $ (i.e.\ where we drop the metrization data) together with the
  directions of all its edges.

  The space of all $n$-marked rational tropical stable maps of a given degree $\Sigma$ in
  $\RR^r$ is denoted by $ \calM_{0,n}(\mathbb{R}^r, \Sigma) $. 
\end{definition}

Since $n$ equals the number of zero-entries in $\Sigma$ and thus can be
deduced from $ \Sigma $, we sometimes drop the subscript and write only $
\calM_{0}(\mathbb{R}^r, \Sigma) $. While all $N$ ends come with markings
$x_1,\ldots,x_N$, only the ends with markings $x_1,\ldots,x_n$ are contracted
(i.e.\ have zero direction) and are thus highlighted in the notation.

\begin{remark}[$ \calM_\bgroup 0,n \egroup(\mathbb R^r,\Sigma) $ as a tropical
    variety]\label{rem-m0nrrdelta}
  We assume $ n \ge 1 $. Then by
  \cite[proposition 4.7]{GKM07}, $ \calM_{0,n}(\mathbb{R}^r,\Sigma) $ is a
  tropical variety, identified with $\calM_{0,N} \times \RR^r$ via the
  map
    \[ \calM_{0,n}(\mathbb{R}^r,\Sigma) \to \calM_{0,N} \times \RR^r, \quad
       (\Gamma,x_1,\dots,x_N,h) \mapsto ((\Gamma,x_1,\dots,x_N), h(x_1)) \]
  which forgets $h$, but records the image $h(x_1)$ of a root vertex. It thus
  inherits the fan structure of $\calM_{0,N}$. In particular, it can be
  embedded via this map into $  \RR^{\binom{N}{2}}/U_N
  \times \RR^r$. When we work with an element of $ \calM_{0,n}(\mathbb{R}^r,
  \Sigma) $ in coordinates, we usually give its coordinates in $
  \RR^{\binom{N}{2}} \times \RR^r$, i.e.\ its image under the distance map and
  the position of the root vertex. If $n=0$ it is still possible to find
  suitable coordinates for $ \calM_{0,n}(\mathbb{R}^r, \Sigma)$ as $\calM_{0,N}
  \times \RR^r$, not by evaluating a marked end but by evaluating for example a
  barycenter \cite[construction 1.2.21]{O13}.
\end{remark}

For each $i=1,\ldots,n$, we have the \emph{evaluation map} 
  \[ \ev_i : \calM_{0,n}(\mathbb{R}^r,\Sigma) \rightarrow \RR^r \]
assigning to a tropical stable map $(\Gamma,x_1,\dots,x_n,h)$ the position
$h(x_i)$ of its $i$-th marked end. It is shown in \cite[proposition 4.8]{GKM07}
that these maps are morphisms of tropical fans.

As above, we will also allow curves in $\RR^r$ where some of the non-contracted
ends are bounded, and write the corresponding moduli spaces as $ \calM_{0,n}'
(\mathbb{R}^r,\Sigma) $.

In the following, we will compute several intersection products in cells of
tropical moduli spaces. Since we are often interested in a local situation, we
can restrict to curves of a given combinatorial type $\alpha$. Local
coordinates for the cell of curves of type $\alpha$ are given by the
coordinates of the root vertex and the lengths of each bounded edge. The map
sending a unit vector in these local coordinates to a vector $v_I$ as above is
a unimodular transformation to the vector space spanned by the corresponding
cell in the moduli space. Therefore we can compute lattice indices also in
these local coordinates.

%% file: gluing.tex
\section{The polyhedral complex
  \texorpdfstring{$\calM_{0,n}(L,\Sigma)$}{M 0,n(L,Sigma)}
  and its gluing weights} \label{sec-m0nldelta}

For the whole section, let $L\subset \R^r$ be a smooth tropical curve as in
example \ref{ex-curves}, and let $ \Sigma $ be the degree of a tropical
$n$-marked stable map to $ \R^r $. We want to define a moduli space $
\calM_{0,n}(L,\Sigma) $ of tropical $n$-marked stable maps to $L$ as a tropical
variety. Let us first construct this space as a polyhedral complex, and then
define its weights in the next subsection.

\subsection{The polyhedral complex
  \texorpdfstring{$\calM_{0,n}(L,\Sigma)$}{M 0,n(L,Sigma)}}
  \label{subsec-ascomplex}

We have already mentioned that not all stable maps with image in $L$ will be
allowed in $\calM_{0,n}(L,\Sigma)$. Instead, we have to impose the so-called
Riemann-Hurwitz condition that we introduce now. As we will see in
construction \ref{con-classicalspace}, it corresponds to a local realizability
condition.

\begin{notation}[Covering degrees] \label{not-local}
  Let $(\Gamma,x_1,\dots,x_N,h)\in \calM_{0,n}(\R^r,\Sigma)$ satisfy
  $h(\Gamma) \subset L$ as sets. As $L$ is irreducible we have $h_\ast(\Gamma)
  = d \cdot L$ for some integer $d$ (which depends only on $\Sigma$). We call
  $d$ the \emph{covering degree} of the stable map and denote it by
  $\deg(\Sigma)$.

  For a vertex $V$ of $ \Gamma $, the \emph{local degree} $ \Sigma_V $ at $V$
  is the collection of the directions of its adjacent edges, labeled in an
  arbitrary way starting with the zero directions. We let $N_V=|\Sigma_V|$ and
  $n_V$ the number of zero directions in $\Sigma_V$ (which may come from marked
  ends or contracted bounded edges). The local covering degree will be denoted
  $ d_V = \deg(\Sigma_V) $.
\end{notation}

\begin{definition}[Riemann-Hurwitz number] \label{def-RH}
  Let $(\Gamma, x_1,\ldots,x_N,h)\in \calM_{0,n}(\R^r,\Sigma)$ satisfy
  $h(\Gamma) \subset L$. We define the \emph{Riemann-Hurwitz number} of a
  vertex $V$ of $ \Gamma $ with image $ W = h(V) $ as 
    $$ \mbox{RH}(V) = N_V-n_V-d_V\cdot (\val(W)-2)-2 $$
  (where $ \val W = 2 $ if $W$ lies in the interior of an edge of $L$).
  Note that it depends only on the combinatorial type of the stable map.
\end{definition}

The Riemann-Hurwitz number gives a realizability condition for tropical stable
maps to smooth curves. It appears e.g.\ in
\cite[definition 2.2]{BBM10}, \cite[proposition 2.4]{Cap12},
\cite[section 3.2.2]{CMR14}, and \cite[definition 3.11]{BruMa}.

\begin{definition}[$\calM_\bgroup 0,n \egroup(L,\Sigma) $ as a polyhedral
    complex] \label{def-ml}
  Let $\alpha$ be a combinatorial type of tropical stable maps in
  $\calM_{0,n}(\R^r,\Sigma)$. We denote the subset of maps
  $(\Gamma,x_1\ldots,x_N,h)$ of type $\alpha$ and satisfying $ h(\Gamma)
  \subset L $ by $\calM(\alpha)$; this is easily seen to be a partially open
  polyhedron. Let $\calM_{0,n}(L,\Sigma)$ be the set of all such cells
  $ \calM(\alpha) $ with $\mbox{RH}(V)\geq 0$ for all vertices $V$ in $ \alpha
  $; this is a polyhedral complex \cite{BM13}.
\end{definition}

Note that this definition of $\calM_{0,n}(L,\Sigma)$ formally differs from the
one used in \cite{GO14} in order to make it compatible with the literature
mentioned above. In \cite[definition \refx{modulx:def-mx}]{GO14}, more cells
are included a priori, but they obtain weight zero in the gluing construction
of subsection \ref{subsec-gluing}. 

\begin{remark}[Dimension of $ \calM_\bgroup 0,n \egroup(L,\Sigma) $]
    \label{rem-ml-dim}
  By an easy generalization of \cite[lemma 2.14]{BM13}, it follows that
  $\calM_{0,n}(L,\Sigma)$ is pure of dimension $|\Sigma|- \deg(\Sigma)\cdot
  \sum_{W \in L} (\val(W)-2) -2$. The maximal cells correspond to combinatorial
  types such that
  \begin{itemize}
  \item each vertex mapping to a vertex of $L$ satisfies $\mbox{RH}(V)=0$,
  \item each vertex mapping to an edge of $L$ is $3$-valent, and
  \item no edge is contracted to a vertex.
  \end{itemize}
\end{remark}

More precisely, we have:

\begin{lemma}[Dimension of cells of $\calM_\bgroup 0,n \egroup(L,\Sigma) $]
    \label{lem-dim}
  Let $\alpha$ be a combinatorial type in $\calM_{0,n}(L,\Sigma)$. The
  dimension of the corresponding cell $\calM(\alpha)$ equals the number of
  vertices mapping to edges of $L$ plus the number of bounded edges mapping to
  vertices of $L$.
\end{lemma}

Intuitively, this holds true since we can independently vary the length of each
bounded edge mapping to a vertex without leaving the cell of a combinatorial
type, as well as the lengths of edges adjacent to a vertex mapping to an edge,
in the appropriate way that ``moves'' the vertex along the edge.

\subsection{The gluing construction for moduli spaces}\label{subsec-gluing}

In this section, we want to equip $\calM_{0,n}(L,\Sigma)$ with weights
satisfying the balancing condition, to make it a tropical variety. To do this,
we review the general technique developed in \cite{GO14}, adapted to the case
when the target of the stable maps is a smooth curve. The idea is to construct
the tropical moduli spaces by a gluing procedure from local moduli spaces for
the vertices. This construction depends on a condition: all vertices appearing
in a combinatorial type of the moduli space are required to be ``good''. We
start by repeating the relevant definitions in the case of smooth curves. 

\begin{notation}[Links of vertices] \label{not-lv}
  Let $(\Gamma,x_1,\ldots,x_n,h)\in \calM_{0,n}(L,\Sigma) $, and let $V$ be a
  vertex of $ \Gamma $. We denote by $L_V$ the link of $L$ around $h(V)$.
  Generalizing the notation of example \ref{ex-curves}, we denote a point by
  $L_0^0$, so that $L_V$ is (an affine shift of a unimodular transformation of)
  $L_r^q\times \R^s$, where $r+s=1$ and $q=0$ if $r=0$. Hence we have $
  (r,s)=(1,0) $ if $V$ maps to a vertex of $L$ (of valence $ q+1 $), and $
  (r,s)=(0,1) $ if $V$ maps to an edge. Note that there is an associated local
  moduli space $ \calM_0(L_V,\Sigma_V) $.
\end{notation}

\begin{definition}[Resolution dimension] \label{def-rdim}
  For a tropical stable map $(\Gamma,x_1,\ldots,x_n,h)\in \calM_{0,n}(L,\Sigma) $, let
  $V$ be a vertex of $\Gamma$ with image $ W = h(V) \in L$. As in notation
  \ref{not-lv}, we have $ L_V\cong L_r^q\times \R^s$ with $r+s=1$ and $q=0$ if
  $r=0$. Treating again a point on an edge of $L$ as a $2$-valent vertex, we
  define the \emph{resolution dimension} of $V$ as
    $$ \rdim(V)= N_V-d_V\cdot (\val(W)-2)+r-3 $$
  and the \emph{classification number} as
    $$ c_V =  N_V+r \quad \in \N. $$
\end{definition}

\begin{remark}[Dimension of local moduli spaces] \label{rem-local-dim}
  By the dimension formula, we see that the local moduli
  space at $V$ has dimension $ \dim \calM_0(L_V,\Sigma_V) = \rdim(V) + s $,
  where again $ L_V \cong L^q_r \times \RR^s $. As this moduli space has an
  $s$-dimensional lineality space coming from shifting the curves along $ \RR^s
  $, the resolution dimension of $V$ is just the dimension of the local moduli
  space at $V$ modulo its lineality space.
\end{remark}

\begin{remark}[Dimension of $ \calM_\bgroup 0,n \egroup (L,\Sigma) $ in terms
    of resolution dimensions] \label{rem-rdim-dim}
  Let $ \alpha $ be a combinatorial type in $ \calM_{0,n}(L,\Sigma) $, and
  assume that $ \alpha $ has $s$ vertices mapping to an edge in $L$ (i.e.\ so
  that the corresponding link is $ L^0_0 \times \RR $). Adding up the
  resolution dimensions of all vertices in $ \alpha $, we obtain by remark
  \ref{rem-ml-dim}
    \[ \sum_V \rdim(V) + s = \dim \calM_{0,n}(L,\Sigma). \]
\end{remark}

\begin{remark}
  Note that $ \rdim(V) $ and $ \RH(V) $ are very similar: in fact, $ \rdim(V) $
  is just $ \RH(V) $ with additional contributions
  \begin{enumerate}
  \item $ n_V $ of the number of contracted edges at $V$, and
  \item $ -1 $ if $V$ maps to an edge of $L$.
  \end{enumerate}
  In particular, the condition $ \RH(V) \ge 0 $ of definition \ref{def-ml}
  also implies $ \rdim(V) \ge 0 $ (otherwise we would have $ \RH(V)=0 $ and $
  \rdim(V)=-1 $, i.e.\ $V$ maps to an edge, $ N_V=2 $, and $ n_V=0 $, which is
  a contradiction since we do not allow $2$-valent vertices).

  The reason to introduce the numbers of definition \ref{def-rdim} is that they
  are used in the recursive definition of good vertices and the weights of $
  \calM_{0,n}(L,\Sigma) $ below. For this construction we start with the case
  of resolution dimension $0$ and pass to the general case by gluing. The
  initial case is obtained by passing to the corresponding situation in
  algebraic geometry and considering (algebraic) Hurwitz numbers.
\end{remark}

\begin{construction}[Algebraic moduli spaces for a vertex]
    \label{con-classicalspace}
  Let $V$ be a vertex of a combinatorial type in $ \calM_{0,n}(L,\Sigma) $ such
  that $ L_V\cong L^q_1 $. Up to unimodular transformation, $ \Sigma_V =
  (\delta_1,\dots,\delta_{N_V})$ is a degree of tropical stable maps to $\R^q$ with ends in
  the directions of $L^q_1$. We decompose $\{1,\dots,N_V\}$ into a partition
  $\eta_0,\dots,\eta_q$ and $\eta$, where
    $$ \eta_i=\{j|\delta_j= -m_j e_i\ \mbox{ for some }
       m_j\in\N_{>0} \} $$
  and $\eta=\{j|\delta_j=0\}$. This also uniquely defines the values
  $m_j$ as the weights of the edges adjacent to $V$.

  To construct an algebraic moduli space for $V$, fix $ q+1 $ distinct
  points $ P_0,\dots,P_q $ on the complex projective line $ \PP^1 $.
  Inside the well-known moduli stack $ \overline{M}_{0,N_V}(\mathbb{P}^1,d_V) $
  of $ N_V $-marked degree-$d_V$ rational stable maps to $ \PP^1 $, consider
  the substack $ M(\Sigma_V) $ of all smooth stable maps $ \mathcal{C} =
  (C,x_1,\dots,x_{N_V},\pi) $ such that $ \pi^* P_i = \sum_{j\in\eta_i}
  m_j x_j $ for all $ i=0,\dots,q $, i.e.\ such that the ramification
  profile of $ \pi $ over $ P_0,\dots,P_q $ is as specified by $ \Sigma_V $. We
  denote its closure inside $\overline{M}_{0,N_V}(\mathbb{P}^1,d_V)$ by $
  \overline{M}(\Sigma_V) $, and its boundary by $\partial M(\Sigma_V) =
  \overline{M}(\Sigma_V) \setminus M(\Sigma_V)$. Its dimension is
    \[ \dim M(\Sigma_V) = 2 d_V - 2 + N_V  - d_V \cdot (q+1) = \rdim(V). \]
\end{construction}

\begin{construction}[The case of resolution dimension $0$] \label{con-weights}
  Let $V$ be a vertex of a combinatorial type in $\calM_{0,n}(L,\Sigma)$ with $
  \rdim(V) = 0 $, where $ L_V \cong L^q_r \times \R^s $ as above. Then
  $ \dim \calM_0(L_V,\Sigma_V) \cong \RR^s $ by remark \ref{rem-local-dim},
  i.e.\ the local moduli space at $V$ consists of only one cell. We make it
  into a tropical variety by giving it the following local weight $ \omega_V $,
  depending on whether $V$ maps to a vertex or to an edge of $L$.
  \begin{enumerate} \parindent 0mm \parskip 1ex plus 0.3ex
  \item If $L_V\cong L^q_1$, the algebraic moduli space $ \overline M(\Sigma_V)
    $ of construction \ref{con-classicalspace} has dimension zero. We define
    the \emph{local weight} of $V$ to be $ \omega_V := \deg \overline
    M(\Sigma_V) $; i.e.\ the number of points in $ \overline M(\Sigma_V) $,
    counted with weight $|\Aut(\pi)|^{-1}$ as we work with a stack. This number
    is also called the \emph{(marked) Hurwitz number} and denoted $ H(\Sigma_V)
    $.
  \item If $ L_V\cong L^0_0\times \R $, the dimension condition implies $ N_V=3
    $. In this case, we set $ \omega_V := 1 $.
  \end{enumerate}
\end{construction}

In fact, the second case could be treated similarly to the first one by
introducing a rubber variant of the moduli space $ \overline M(\Sigma_V) $. We
avoid this formulation for the sake of simplicity.

Let us now describe the gluing construction that gives the local moduli space $
\calM_0(L_V,\Sigma_V) $ of a vertex $V$ the structure of a tropical variety if
$ \rdim(V)>0 $. In the following, any combinatorial type occurring in $
\calM_0(L_V,\Sigma_V) $ will be called a \emph{resolution} of $V$. For a
combinatorial type $ \alpha $ occurring in a moduli space we denote by $
\calN(\alpha) $ the ``neighborhood of $ \alpha $'', i.e.\ the union of all
cells $ \calM(\beta) $ whose closure intersects $ \calM(\alpha) $.

Definition \ref{def:good} of a good vertex and the following gluing
construction \ref{con:gluing} depend on each other and work in a combined
recursion on the classification number of vertices. The following definition of
a good vertex thus assumes that good vertices of lower classification number
are already defined recursively. Moreover, for every combinatorial type
$\alpha$ in a local moduli space $\calM_0(L_V,\Sigma_V) $ all of whose vertices
have smaller classification number and are good it assumes that there is a
gluing cycle in the neighborhood $\calN(\alpha)$ from construction
\ref{con:gluing}.

\begin{definition}[{Good vertices \cite[definition \refx{modulx:def:good}]{GO14}}]
  \label{def:good}

Let $V$ be a vertex of a (local) tropical stable map in $\calM_0(L_V,\Sigma_V) $, so that in particular $\rdim(V)\geq 0$ . The vertex $V$ is called \textit{good} if the following holds:
  \begin{enumerate}
  \item Every vertex of every resolution $\alpha$ of $V$ in  $\calM_0(L_V,\Sigma_V) $ (which has classification number smaller than $c_V$ by
    \cite[lemma \refx{modulx:lem-classnumber}]{GO14}) is good (so that a gluing cycle is defined on $\calN(\alpha)$ by construction \ref{con:gluing}).
\item If $\rdim(V)>0$ the maximal types in $\calM_0(L_V,\Sigma_V) $ are resolutions of $V$. We let $\calM_0(L_V,\Sigma_V) $ be a weighted polyhedral complex by defining the weights on maximal cells $\calM(\alpha)=\calN(\alpha)$ using the gluing construction \ref{con:gluing}. If $\rdim(V)=0$, we equip the unique cell of $\calM_0(L_V,\Sigma_V) $ with the weight of construction \ref{con-weights}. We require that the space $\calM_0(L_V,\Sigma_V) $ is a tropical cycle with these weights.
\item For every resolution $\alpha$ of $V$ in $\calM_0(L_V,\Sigma_V) $ and
  every maximal type $\beta$ such that $ \overline{\calM(\beta)} $ contains $\calM(\alpha)$ in $\calM_0(L_V,\Sigma_V) $ ($\beta$ is then also a resolution of $V$), the weight of $\beta$ is the same in the gluing cycles $\calN(\alpha)$ and $\calN(\beta)$.
  \end{enumerate}
\end{definition}

In the following review of the gluing construction from
\cite[construction \refx{modulx:con:gluing}]{GO14}, we omit some of the
technical details for the sake of clarity.

\begin{construction}[The gluing construction for a combinatorial type $ \alpha
    $] \label{con:gluing}
  Fix a (not necessarily maximal) combinatorial type $\alpha$ of curves in
  $\calM_{0,n}(L,\Sigma)$ and assume that all its vertices are good. We will
  construct weights on the maximal cells of the neighborhood $ \calN(\alpha) $
  such that this partially open polyhedral complex becomes a tropical cycle.
  In particular, if $ \alpha $ is already maximal this defines a weight on $
  \calM(\alpha) = \calN(\alpha) $.

  We cut each bounded edge of $\alpha$ at some point in its interior, and in
  addition introduce lengths for all ends. This yields a set of connected
  components $ \alpha_V $, each containing only one vertex $V$, edges of
  directions $ \Sigma_V $, and (now bounded) ends labeled by an index set $ I_V
  $.

  For every such vertex $V$, consider the local moduli space $
  \calM_{0}(L_V,\Sigma_V)$, which is a tropical variety since $V$ is good. We
  introduce lengths on all ends of $ \Sigma_V $, obtaining a moduli space $
  \calM'_{0}(L_V,\Sigma_V) $ (of which $\calM_{0}(L_V,\Sigma_V)$ is a
  quotient) as in subsection \ref{subsec-m0nrd}. Each bounded end $i\in I_V$ is
  mapped to an edge or vertex of $L$ that we denote by $\sigma_i$. We consider
  the open subcomplex of $\calM'_{0}(L_V,\Sigma_V)$ of all curves for which the
  evaluation at $i$ still lies in $ \sigma_i $, i.e.\ the partially open
  tropical subvariety
    $$ \calM_V := \bigcap_{i\in I_V}\ev_i^{-1}(\sigma_i) $$
  of $\calM'_{0}(L_V,\Sigma_V)$.

  Now we want to glue these pieces $\calM_V$ back together. Consider a bounded
  edge $e$ of $ \alpha $ adjacent to two vertices $V_1(e)$ and $V_2(e)$, and
  denote the two bounded ends produced by cutting $e$ by $i_1(e) \in I_{V_1(e)}
  $ and $i_2(e)\in I_{V_2(e)}$, where $ \sigma_{i_1(e)} = \sigma_{i_2(e)} =:
  \sigma_e $. There is a corresponding evaluation map
    $$ \ev_e := (\ev_{i_1(e)} \times \ev_{i_2(e)}) :
       \prod_V \calM_V \longrightarrow \sigma_e \times \sigma_e $$
  at the endpoints of these two bounded ends in the factors for $ V_1 $ and $
  V_2 $. To impose the condition that these ends fit together to form the edge
  $e$ we need to pull back the diagonal $\Delta_{\sigma_e}$ via $\ev_e$
  \cite[appendix \refx{modulx:sec:diagonal}]{GO14}. We abbreviate all these
  pull-backs by
    $$ \ev^*(\Delta_{L}) \cdot \prod_V \calM_V :=
       \prod_e \ev_e^*\Delta_{\sigma_e} \cdot
       \prod_V \calM_V, $$
  where $e$ runs over all bounded edges $e$ of $\alpha$. By construction, this
  cycle consists of stable map pieces that glue back to a stable map in $
  \calM_{0,n}(L,\Sigma) $. However, it also carries the superfluous
  information on the position of the gluing points. To get rid of this we
  apply the quotient map $q$ by the lineality space generated by the
  appropriate differences of vectors taking care of the lengths of the bounded
  ends, and by the vectors taking care of ends which should be unbounded. We
  finally use a morphism $f$ identifying a stable map glued from pieces
  with the corresponding element in $\calM_{0,n}(L,\Sigma)$, where we use the
  distance and barycentric coordinates mentioned in remark
  \ref{rem-m0nrrdelta}. Hence we get a partially open tropical cycle
    $$ f_*q\left[\ev^*(\Delta_{L})\cdot\prod_V\calM_V\right]
       \quad\text{in}\quad
       \calM_{0,n}(L,\Sigma). $$
  Its weights on the maximal cells of $ \calM_{0,n}(L,\Sigma) $ will be called
  the \textit{gluing weights}. It is easy to see that the gluing morphism $f$
  is unimodular and induces a bijection of cells. In particular, the weight of
  a maximal cell in $f_*q\left[\ev^*(\Delta_{L})\cdot\prod_V\calM_V\right]$
  is equal to the weight of $ \ev^*(\Delta_{L})\cdot\prod_V\calM_V $ in the
  corresponding cell of $ \prod_V \calM_V $. By remark \ref{rem-latticeindex},
  it can be computed as the greatest common divisor of the maximal minors of a
  matrix whose rows represent the differences $ \ev_{i_1(e)} - \ev_{i_2(e)} $
  in local coordinates.
\end{construction}

\begin{example} \label{ex-evl2}
  Let $L=L^2_1$ be a tropical line in $\R^2$ and let $\alpha$ be a
  combinatorial type of degree-$\Sigma$ curves in $L^2_1$ as shown below on the
  left (where the directions of the edges indicate their images in $ \R^2 $).
  Then $ \rdim(V_0)=0 $. We assume in addition that $\rdim(V_1)=0$.

  \begin {center} \input {pics/mult1} \end {center}

  We cut the unique bounded edge $e$ of weight $ d_1 $, obtaining two bounded
  ends that we denote $f$ and $f'$. By the assumption on the resolution
  dimension, the local moduli spaces for $V_0$ and $V_1$ consist of only one
  cell each, and we can explicitly describe isomorphisms to open polyhedra in
  some $\R^k$ as follows. The space $\calM_{V_0}$ is isomorphic to $\R_{>0}^{2}
  $, where one coordinate that we denote by $l_{f}$ corresponds to the length
  of the bounded end, and the other that we call $x_{V_0}$ to the position of
  the image of $V_0$ on the corresponding ray of $L$. The space $\calM_{V_1}$
  is $\R_{>0}$ with coordinate $l_{f'}$ corresponding to the length of its
  bounded end. By construction \ref{con-weights}, the weight of $\calM_{V_1}$
  is the Hurwitz number $ \omega_{V_1} = H(\Sigma_{V_1}) $, whereas
  $\calM_{V_0}$ has weight $1$. Using these coordinates, we can pull back the
  diagonal of $L$ as $ \ev_e^* \max\{x-y,0\}=\max\{\ev_{f}-\ev_{f'},0\}$, where
  $x,y$ are the coordinates of $ L^2 $ on the left ray. By remark
  \ref{rem-latticeindex}, the weight of $\ev^*\Delta_{L} \cdot (\calM_{V_0}
  \times\calM_{V_1})$ equals the weight of $\calM_{V_0}\times\calM_{V_1}$ times
  the greatest common divisor of the maximal minors of the matrix
    \[ \begin {array}{l|cccc}
       & x_{V_0} & l_{f} & l_{f'}   \\ \hline
       \ev_{f}-\ev_{f'} &
         1 & -d_1 & -d_1,   \\
       \end {array} \]
  which is $1$. Hence the cell corresponding to $\alpha$ in $\calM_0(L,\Sigma)
  $ has weight $\omega_\alpha = H(\Sigma_{V_1}) $. The analogous result holds for $
  L=L^q_1 $ for all $q$.
\end{example}

\begin{example} \label{ex-evl}
  Let $L=L^2_1$ be a tropical line in $\R^2$ again, and let $\alpha$ be the
  combinatorial type of degree-$ \Sigma $ curves mapping to $L^2_1$ depicted
  below, with $ V_1 $ and $ V_2 $ mapping to the vertex of $L$. As above, we
  then have $\rdim(V_0)=0$, and assume in addition that $\rdim(V_1)=\rdim(V_2)
  = 0$.

  \begin {center} \input {pics/mult2} \end {center}

  We cut the two edges $e_1$ of weight $d_1$ and $e_2$ of weight $d_2$,
  obtaining four new bounded ends that we denote by $f_i$ and $f_i'$ for
  $i=1,2$. As before, each local moduli space consists of only one
  cell. The space $\calM_{V_0}$ is isomorphic to $\R_{>0}^{3}$, where two
  coordinates ($l_{f_1}$ and $l_{f_2}$) correspond to the lengths of the
  bounded ends and one ($x_{V_0}$) to the position of the image of $V_0$ on the
  corresponding ray of $L$. By construction \ref{con-weights}, it is equipped
  with weight $\omega_{V_0}=1$. Similarly, $\calM_{V_i}$ for $ i=1,2 $ is
  isomorphic to $\R_{>0}$, where the coordinate $l_{f'_i}$ is given by the
  length of the bounded end, and equipped with the appropriate Hurwitz number $
  \omega_{V_i}=H(\Sigma_{V_i})$ as weight. As in the previous example, pulling
  back the diagonal of $L^2$ twice and using remark \ref{rem-latticeindex}, we
  deduce that the weight of $ \ev^*\Delta_{L} \cdot (\calM_{V_0}\times
  \calM_{V_1} \times \calM_{V_2}) $ equals the weight of $ \calM_{V_0}\times
  \calM_{V_1} \times \calM_{V_2} $ times the greatest common divisor of the
  maximal minors of the matrix
    \[ \begin {array}{l|ccccc}
         & x_{V_0} & l_{f_1} & l_{f_2} & l_{f_1'} & l_{f'_2}  \\ \hline
         \mbox {$\ev_{f_1}-\ev_{f'_1}$} &
           1 & -d_1 & 0 & -d_1 & 0  \\
         \mbox {$\ev_{f_2}-\ev_{f'_2}$} &
           1 & 0 & -d_2 & 0 & -d_2, \\
       \end {array} \]
  which is $ \gcd(d_1,d_2) $. Thus the weight of the cell corresponding to
  $\alpha$ in $\calM_0(L,\Sigma)$ equals
    $$ \omega_\alpha = \gcd(d_1,d_2) \, \omega_{V_0} \omega_{V_1}
       \omega_{V_2} = \gcd(d_1,d_2) \cdot H(\Sigma_{V_1}) \cdot
       H(\Sigma_{V_2}). $$
  As in example \ref{ex-evl2}, the same result holds for $ L = L^q_1 $ for all
  $q$.
\end{example}

We end this section by stating the main result of \cite{GO14}, together with a
lemma that provides a major simplification for checking the requirements of the
following theorem:

\begin{theorem}[{The gluing theorem
    \cite[corollary \refx{modulx:cor-balancing}]{GO14}}] \label{thm-glue}
  Assume that all vertices $V$ that can possibly occur in combinatorial types
  of the moduli space $\mathcal{M}_{0,n}(L,\Sigma)$ are good. Then the
  gluing construction is well-defined for all these combinatorial types. In
  particular, $\calM_{0,n}(L,\Sigma)$ is a tropical variety.
\end{theorem}

\begin{lemma}[{Restriction to resolution dimension one
    \cite[corollary \refx{modulx:cor:good1dim}]{GO14}}]
    \label{lem-good1dim}
  If all vertices $V$ of combinatorial types of $\mathcal{M}_{0,n}(L,\Sigma)$
  with $\rdim(V)=1$ are good, then all vertices are good. 
\end{lemma}

%% file: pics/mult1.tex
\begin{picture}(0,0)%
\includegraphics{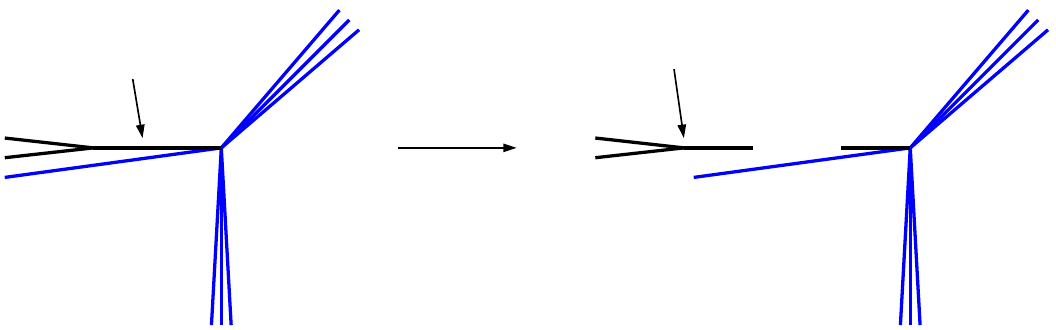}%
\end{picture}%
\setlength{\unitlength}{4144sp}%
\begingroup\makeatletter\ifx\SetFigFont\undefined%
\gdef\SetFigFont#1#2#3#4#5{%
  \reset@font\fontsize{#1}{#2pt}%
  \fontfamily{#3}\fontseries{#4}\fontshape{#5}%
  \selectfont}%
\fi\endgroup%
\begin{picture}(4814,1507)(339,-893)
\put(766,-219){\makebox(0,0)[b]{\smash{{\SetFigFont{10}{12.0}{\familydefault}{\mddefault}{\updefault}{\color[rgb]{0,0,0}$V_0$}%
}}}}
\put(946,321){\makebox(0,0)[b]{\smash{{\SetFigFont{10}{12.0}{\familydefault}{\mddefault}{\updefault}{\color[rgb]{0,0,0}weight $d_1$}%
}}}}
\put(946,479){\makebox(0,0)[b]{\smash{{\SetFigFont{10}{12.0}{\familydefault}{\mddefault}{\updefault}{\color[rgb]{0,0,0}edge $e$}%
}}}}
\put(4546,-196){\makebox(0,0)[lb]{\smash{{\SetFigFont{10}{12.0}{\familydefault}{\mddefault}{\updefault}{\color[rgb]{0,0,1}$V_1$}%
}}}}
\put(3421,366){\makebox(0,0)[b]{\smash{{\SetFigFont{10}{12.0}{\familydefault}{\mddefault}{\updefault}{\color[rgb]{0,0,0}$V_0$}%
}}}}
\put(3646,  6){\makebox(0,0)[b]{\smash{{\SetFigFont{10}{12.0}{\familydefault}{\mddefault}{\updefault}{\color[rgb]{0,0,0}$f$}%
}}}}
\put(4366,  6){\makebox(0,0)[b]{\smash{{\SetFigFont{10}{12.0}{\familydefault}{\mddefault}{\updefault}{\color[rgb]{0,0,0}$f'$}%
}}}}
\put(1396,-196){\makebox(0,0)[lb]{\smash{{\SetFigFont{10}{12.0}{\familydefault}{\mddefault}{\updefault}{\color[rgb]{0,0,1}$V_1$}%
}}}}
\end{picture}%

%% file: pics/mult2.tex
\begin{picture}(0,0)%
\includegraphics{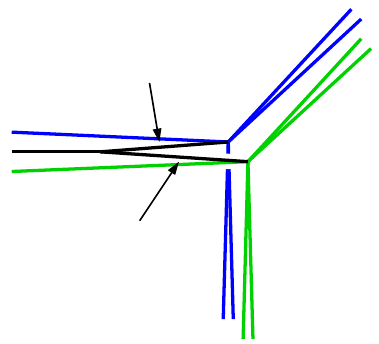}%
\end{picture}%
\setlength{\unitlength}{4144sp}%
\begingroup\makeatletter\ifx\SetFigFont\undefined%
\gdef\SetFigFont#1#2#3#4#5{%
  \reset@font\fontsize{#1}{#2pt}%
  \fontfamily{#3}\fontseries{#4}\fontshape{#5}%
  \selectfont}%
\fi\endgroup%
\begin{picture}(1719,1564)(2782,-983)
\put(3826,  6){\makebox(0,0)[rb]{\smash{{\SetFigFont{10}{12.0}{\familydefault}{\mddefault}{\updefault}{\color[rgb]{0,0,1}$V_1$}%
}}}}
\put(3961,-286){\makebox(0,0)[lb]{\smash{{\SetFigFont{10}{12.0}{\familydefault}{\mddefault}{\updefault}{\color[rgb]{0,.82,0}$V_2$}%
}}}}
\put(3241,-264){\makebox(0,0)[b]{\smash{{\SetFigFont{10}{12.0}{\familydefault}{\mddefault}{\updefault}{\color[rgb]{0,0,0}$V_0$}%
}}}}
\put(3466,276){\makebox(0,0)[b]{\smash{{\SetFigFont{10}{12.0}{\familydefault}{\mddefault}{\updefault}{\color[rgb]{0,0,0}weight $d_1$}%
}}}}
\put(3466,434){\makebox(0,0)[b]{\smash{{\SetFigFont{10}{12.0}{\familydefault}{\mddefault}{\updefault}{\color[rgb]{0,0,0}edge $e_1$}%
}}}}
\put(3421,-736){\makebox(0,0)[b]{\smash{{\SetFigFont{10}{12.0}{\familydefault}{\mddefault}{\updefault}{\color[rgb]{0,0,0}weight $ d_2 $}%
}}}}
\put(3421,-579){\makebox(0,0)[b]{\smash{{\SetFigFont{10}{12.0}{\familydefault}{\mddefault}{\updefault}{\color[rgb]{0,0,0}edge $ e_2 $}%
}}}}
\end{picture}%

%% file: curves.tex
\section{One-dimensional moduli spaces of rational covers of smooth tropical
  curves} \label{sec-1dim}

Throughout this section, let $V$ be a vertex of a combinatorial type in
$\calM_{0,n}(L,\Sigma)$ with $ \rdim(V)=1 $. Our aim is to show that $V$ is
good, so that we can apply lemma \ref{lem-good1dim} and the gluing theorem
\ref{thm-glue} to deduce theorem \ref{thm1}. We continue to use the notation of
section \ref{sec-m0nldelta}. Moreover, let $ I_V $ be the set of labels of the
ends in the local moduli space $ \calM_{0}(L_V,\Sigma_V) $, so that $ \calM_{0}
(L_V,\Sigma_V) = \calM_{0,I_V}(L_V,\Sigma_V) $. As in construction
\ref{con-classicalspace}, let $ m_j \in \N_{>0} $ be the weight of the end
$ j \in I_V $.

To prove that $V$ is good, we have to show by definition \ref{def:good} that
\begin{enumerate}
\item[(1)] every vertex appearing in a non-trivial resolution in $
  \calM_{0}(L_V,\Sigma_V)$ is good;
\item[(2)] $\calM_{0}(L_V,\Sigma_V)$ is a tropical variety with the gluing
  weights; and
\item[(3)] for every non-trivial resolution $ \alpha $ of $V$, the weight of
  each maximal cell in the neighborhood $ \calN(\alpha) $ is the same no matter
  if we apply the gluing construction for $ \alpha $ or just for this maximal
  cell.
\end{enumerate}
Assume first that $V$ maps to an edge of $L$, so that $ L_V \cong L_0^0\times
\R $. Then $ \rdim(V)=1 $ implies $ N_V=4 $, hence the possible resolutions are
just the usual resolutions of a $4$-valent vertex. Also, any gluing weight is
just $1$, and the balancing condition is satisfied --- this is just the usual
balancing condition of $\calM_{0,4}$. It follows that $V$ is good.

We can thus assume now that $V$ maps to a vertex of $L$, so that $ L_V \cong
L_1^q $. By Remark \ref{rem-local-dim}, this means that $ \dim \calM_{0}
(L_V,\Sigma_V) = 1 $. In particular, every resolution of $V$ corresponds
already to a maximal cell of the local moduli space, which implies that
condition (3) above is trivially satisfied. Moreover, lemma \ref{lem-dim}
implies that every non-trivial resolution of $V$ has at least one vertex
mapping to an edge of $L$, or a bounded edge contracted to a vertex. In the former case,
remark \ref{rem-rdim-dim} then shows that all vertices in this resolution must
have resolution dimension $0$ and are thus good, and the latter case is an
immediate contradiction to remark \ref{rem-ml-dim}. Hence condition (1) is
always satisfied as well, and it only remains to check the balancing condition
(2).

Next, since $ 1 - n_V = \rdim(V)-n_V=\mbox{RH}(V)\geq 0$, we can either have
$n_V=1$ and $\mbox{RH}(V)=0$, or $n_V=0$ and $\mbox{RH}(V)=1$. In the first
case, there is one contracted end, say with the marking $1$, adjacent to the
vertex. In the possible resolutions, this contracted end is adjacent to any
other of the non-contracted ends, leading to a generating vector of the form
$v_{\{1,i\}}$ for the corresponding ray in $\calM_{0}(L_V,\Sigma_V)$. As in
example \ref{ex-evl2}, we can see that any gluing weight equals
$H(\Sigma_V\setminus\{0\})$. We have $\sum_{i=2}^{N_V} v_{\{1,i\}}=0$ in
$\calM_{0}(\R^q,\Sigma_V)$, and hence the balancing condition is satisfied in
this case.

So the only thing left to be done is to study the remaining case, where we have
a vertex $V$ mapping to a vertex of $L$, without contracted ends and having $
\rdim(V)=1 $, and to prove the balancing condition (2) for the $1$-dimensional
local moduli space $ \calM_{0,n}(L,\Sigma) $ in this situation. We start by
listing the possible resolutions of such a vertex, i.e.\ the maximal cones of $ 
\calM_{0,n}(L,\Sigma) $.

\begin{construction}[Resolutions of a vertex with $\rdim(V)=1$]
    \label{con-codim0types}
  Let $V$ be a vertex of a combinatorial type in $\calM_{0,n}(L,\Sigma)$.
  Assume that $V$ maps to $L_V\cong L^q_1$ and satisfies $\rdim(V)=1$ and
  $n_V=0$.

  As $ \dim \calM_{0}(L_V,\Sigma_V) = 1 $, it follows from remark
  \ref{rem-ml-dim} and lemma \ref{lem-dim} that in each (necessarily maximal)
  resolution of $V$, there is one (necessarily $3$-valent) vertex $ V_0 $
  mapping to an edge of $L^q_1$. This vertex can either join two ends or split
  an end, so that we obtain the following two types of resolutions:
 
  \begin {center} \input {pics/res1} \end {center}

  \begin{enumerate} \parindent 0mm \parskip 1ex plus 0.3ex
  \item[(I)] There is exactly one vertex $ V_1 $ mapping to the vertex of
    $ L_V $. The vertex $ V_0 $ is adjacent to two ends $ i,j \in I_V $ and a
    bounded edge of weight $ d_1 = m_i + m_j $ connecting $ V_0 $ to
    $ V_1 $. The ends in $ I_1 := I_V \backslash \{i,j\} $ are adjacent to $
    V_1 $.

    Such a type exists for all choices of ends $i$ and $j$ of the same (primitive) direction.
  \item[(II)] There are exactly two vertices $ V_1,V_2 $ mapping to the vertex
    of $ L_V $. The vertex $ V_0 $ is adjacent to an end $ i \in I_V $ and two
    bounded edges of weights $ d_1,d_2 $ with $ d_1+d_2 = m_i $ connecting
    $ V_0 $ to $ V_1 $ and $ V_2 $, respectively. The two vertices $ V_1 $ and
    $ V_2 $ are adjacent to ends in $ I_1 $ and $ I_2 $, respectively, where $
    I_1 \cup I_2 \cup \{i\} = I_V $.

    Such a type exists for all choices of $i$ and all partitions of $ I_V
    \backslash \{i\} $ into $ I_1 $ and $ I_2 $ for which there is a stable map with
    the above conditions.
  \end{enumerate}

  With the notations of section \ref{sec-prelim}, these types correspond to
  rays of $ \mathcal{M}_{0}(L_V,\Sigma_V) $ generated by the vectors $
  v_{\{i,j\}} $ for type I and $ d_2v_{I_1}+d_1v_{I_2} $ for type II (where the
  latter does not need to be primitive).

  Let us now consider the corresponding algebraic situation, i.e.\ the
  $1$-dimensional algebraic moduli space $\overline M(\Sigma_V)$ of
  construction \ref{con-classicalspace}. By the Riemann-Hurwitz condition, a
  point in the open part $M(\Sigma_V)$ corresponds to a cover with precisely
  one simple ramification which is not marked, and whose image does not
  coincide with one of the points $ P_0,\dots,P_q $ at which we fixed the
  ramification imposed by $\Sigma_V$. The boundary points correspond to
  degenerate covers that we obtain when the additional branch point runs into
  a point $ P_s $ for $ s \in \{0,\dots,q\} $.

  As deformations of covers are always local around special fibers
  \cite[proposition 1.1]{Vak00}, we see that a cover in $\partial M(\Sigma_V)$
  must have exactly one collapsed component, which then has exactly three
  special points. So we have the following two types for the curves in the
  boundary $ \partial M(\Sigma_V) $, which are exactly dual to the tropical
  picture above (see \cite[proposition 3.12]{BM13} for a related statement):

  \begin {center} \input {pics/res2} \end {center}

  Here, $ C_0 $ is the collapsed component, and $ C_k $ for $ k \in \{1,2\} $
  denotes the at most $2$ non-collapsed irreducible components. In the type I
  case, the map $ \pi|_{C_1} $ has order $ d_1 := m_i+m_j $ at the
  singular point of $C$. In the type II case, the orders $ d_1 $ and $ d_2 $ of
  $ \pi|_{C_1} $ and $ \pi|_{C_2} $ at the singular points of $C$ add up to $
  m_i $.
\end{construction}

To check the balancing condition in the $1$-dimensional fan $
\calM_{0}(L_V,\Sigma_V) $, it suffices by lemma \ref{lem-forget} to consider
the situation after applying the various forgetful maps to $ \calM_{0,4} $. We
will do this first in the algebraic and then in the tropical case.

\begin{lemma}[The pull-back of the forgetful map] \label{lem-boundary}
  Let $ \calC \in \partial M(\Sigma_V) $ be a stable map in the boundary of the
  local moduli space of a vertex $V$ as in construction \ref{con-codim0types}.
  Consider the forgetful map $ \ft_I: \overline{M}(\Sigma_V) \to
  \overline{M}_{0,I} \cong \PP^1 $ for a choice of four-element subset
  $ I = \{i,j,k,l \} \subset I_V $. Then the multiplicity $ \ord_\calC \ft_I^*
  (ij|kl) $ of the pullback of the divisor $ (ij|kl) $ on $ \overline{M}_{0,I}
  $ at $ \calC $ equals
  \begin{enumerate}
  \item $1$ if $ \calC $ is of type I, with $ x_i,x_j \in C_0 $ and $ x_k,x_l
    \in C_1 $ or vice versa;
  \item $d_1$ if $ \calC $ is of type II, with $ x_i \in C_0 $, and $ x_j \in
    C_1 $ and $ x_k,x_l \in C_2 $ or vice versa;
  \item $ d_1+d_2 $ if $ \calC $ is of type II, with $ x_i,x_j \in C_1 $ and $
    x_k,x_l \in C_2 $ or vice versa.
  \end{enumerate}

  \begin {center} \input {pics/res3} \end {center}

  These are all cases in which we have a non-zero multiplicity.
\end{lemma}

\begin{proof}
  Since $ x_i,x_j $ and $ x_k,x_l$ must lie on different components after
  applying the forgetful map and $ \calC $ has at least two and at most three
  components, it is obvious that we can only have the three cases stated in the
  lemma. We want to determine the multiplicity of $ \calC $ in $ \ft_I^*(ij|kl)
  $ for each case. By \cite[proposition 1.1]{Vak00}, we may replace our family
  $ \overline M(\Sigma_V) $ of curves around $ \calC $ by another family $
  \overline M $ of curves \'etale locally isomorphic to the original ones around
  the collapsed component. The following picture illustrates the new curve $
  \calC $ after this replacement in each case; the corresponding families are
  described below.

  \begin {center} \input {pics/res4} \end {center}

  \textbf{Case (a):} Let $M$ be the moduli space of all smooth covers $
  (C,x_i,x_j,x_k,x_l,\pi) $ of $ \PP^1 $ of degree $ d_1 = m_i + m_j $
  satisfying 
    $$ \pi^*0 = m_j x_j + m_i x_i, \;\;
       \pi^*\infty = d_1 x_k, \;\;\text{and}\;\;
       \pi(x_l) = 1. $$
  On the source curve $ C \cong \PP^1 $, we set $x_i=0$, $x_j=\infty$, $x_k=1$,
  and $ x_l=(1:w) $ with $ w\in \C^*\setminus\{1\} $. Then every element in
  $M$ can be written as
    $$ \pi(z_0:z_1)=((z_0-z_1)^{d_1}:\lambda z_0^{m_j}z_1^{m_i}) $$
  for $\lambda\in \C^*$ satisfying $\lambda w^{m_i}=(1-w)^{d_1} $. Thus,
  the $1$-dimensional space $M$ is parametrized by those $(\lambda,w)\in
  \C^*\times(\C^*\setminus\{1\})$ with $\lambda w^{m_i}=(1-w)^{d_1}$. The
  non-marked branch point of $ \pi $ can be computed to be at $P =
  (d_1^{d_1}:(-1)^{m_i} m_j^{m_j} m_i^{m_i} \cdot
  \lambda)$, since the equation $\pi(z_0:z_1)= P$ has a double root at $
  (m_j:-m_i)$.

  Hence, in this family the singular curve $ \calC $ in the picture above
  corresponds to the coordinates $ (\lambda,w)=(0,1)$. After inserting this
  point into the family, we obtain $ \overline M\cong \C^* $ via $(\lambda,w)
  \mapsto w$. The divisor $\ft_I^*(ij|kl)$ is given by the function $w-1$,
  which vanishes to order $1$ at $\mathcal{C} $. As $\mathcal{C}$ has no
  automorphisms due to the marked point $ x_l $, we obtain $\ord_\calC
  \ft_I^*(ij|kl)=1$ as claimed.

  \textbf{Case (b):} Now let $M$ be the space of those smooth covers $
  (C,x_i,x_j,x_k,x_l,\pi) $ of $ \PP^1 $ of degree $ d = m_i $ such that
    $$ \pi^*0=d_1x_j+d_2x_k, \;\;
       \pi^*\infty=dx_i, \;\;\text{and}\;\;
       \pi(x_l)=1 $$
  for fixed $ d_1,d_2 $ with $ d_1+d_2=d $. We set $x_k=1$, $x_i=\infty$,
  $x_j=0$, and $x_l=(1:w)$, where $w\in \C^*\setminus\{1\}$. Then every element
  of $M$ can be written as
    \[ \pi(z_0:z_1)=(\lambda z_0^d:(z_0-z_1)^{d_2}z_1^{d_1}), \]
  where $\lambda\in \C^*$ satisfies $\lambda=(1-w)^{d_1}w^{d_2}$. The
  non-marked branch point of a cover $\pi$ can be computed to be at $
  P=(\lambda\cdot d^d:d_1^{d_1}d_2^{d_2})$, since the equation $\pi(z_0:z_1) =
  P$ has a double root at $(d:d_1)$.

  Again, as in the picture above we want to insert the special fiber $
  \mathcal{C}$ over $(\lambda,w) = (0,1)$ to obtain the space $ \overline M $.
  As before, $\overline M\cong\C^*$ via $(\lambda,w)\mapsto w$, and the divisor
  $ \ft_I^*(ij|kl) $ is given by the function $w-1$, which vanishes to order
  $1$ at $\mathcal{C}$. Since $\mathcal{C}$ has $d_1$ automorphisms on $ C_1 $
  (which is totally ramified over $0$ and $ \infty $), we obtain $\ord_\calC
  \ft_I^* (ij|kl)=d_1$.

  \textbf{Case (c):} In this case, we use the previous computations and the
  WDVV equations. Denote by $x_p$ the marked point of $ \calC $ on the
  collapsed component. We consider the moduli space $\overline M$ which is the
  closure of all smooth $(C,x_i,x_j,x_k,x_l,x_p,\pi)$ of degree $ d=m_p $
  such that
    \[ \pi^*0=d_1x_i+d_2x_l,\;\;
       \pi^*\infty=dx_p,\;\; \text{and} \;\;
       \pi(x_j)=\pi(x_k)=1 \]
  for fixed $ d_1,d_2 $ with $ d_1+d_2 = d $. Again, by the Riemann-Hurwitz
  formula this is a $1$-dimensional space, with one non-marked ramification for
  a smooth curve in $ \overline M $. By letting the additional branch point run
  into $0$, $1$ and $\infty$, we can see that $\partial M$ contains the
  following reducible curves:
  \begin{enumerate}
  \item[(1)] a degree-$d_1$ component with $x_i,x_j$ connected to a
    degree-$d_2$ component with $x_k, x_l$ via a collapsed component over $
    \infty $ with $x_p$ (this is the curve in the picture above);
  \item[(2)] a degree-$d_1$ component with $x_i,x_k$ connected to a
    degree-$d_2$ component with $x_j,x_l$ via a collapsed component over $
    \infty $ with $x_p$;
  \item[(3)] one collapsed component over $0$ with $x_i,x_l$ and one degree-$d$
    component with $x_j,x_k,x_p$;
  \item[(4)] one collapsed component over $1$ with $x_j,x_k$ and one degree-$d$
    component with $x_i,x_l,x_p$.
  \end{enumerate}
  The non-collapsed components in types (1), (2), and (3) are all completely
  ramified over two points. In types (1) and (2), exactly one point with no
  ramification is marked, killing the automorphisms. Hence, for each of these
  types (1) and (2) we have one corresponding boundary point in $\overline M$.
  In type (3), the point $x_j$ fixes the automorphisms, but then we have a
  choice to mark any preimage of $1$ but $x_j$ to be $x_k$. Hence there are $
  d-1 $ boundary points corresponding to a cover of type (3). For type (4), a
  computation of the corresponding Hurwitz number shows that there is a unique
  such cover, so that we have one such boundary point in $\overline M$.

  By the WDVV equations for $\ft_I:M \to \overline{M}_{0,4}$, we have $
  \ft_I^*(ij|kl)= \ft_I^*(il|kj)$. The left side of this equation is obviously
  supported on the boundary point of type (1) that we are interested in,
  whereas the right side is supported on all boundary points of type (3) or
  (4). The multiplicity of $\ft_I^*(il|kj)$ is $1$ at each such boundary point
  by our former computation. As there are $d$ such boundary points, in total we
  obtain $\ord_\calC \ft_I^*(ij|kl)=d$.
\end{proof}

\begin{theorem}[One-dimensional moduli spaces $\calM_0(L_V,\Sigma_V)$]
    \label{thm-local-balancing}
  Let $V$ be a vertex as in construction \ref{con-codim0types}: mapping to $
  L_V\cong L^q_1 $ and satisfying $\rdim(V)=1$ and $n_V=0$. Then $
  \calM_{0}(L_V,\Sigma_V) $ with the weights obtained from the gluing
  construction \ref{con:gluing} is a one-dimensional balanced fan. In
  particular, $V$ is good.
\end{theorem}

\begin{proof}
  The rays of $ \calM_0(L_V,\Sigma_V) $ are given by the combinatorial types
  $ \alpha $ of construction \ref{con-codim0types}. With the notation used
  there, we can take as spanning vectors for these rays $ u_\alpha =
  v_{\{i,j\}} $ in a type I case and $ u_\alpha = d_2 v_{I_1} + d_1 v_{I_2} $
  in a type II case. As the integer length of these vectors is $1$ and $
  \gcd(d_1,d_2) $, respectively, it follows from examples \ref{ex-evl2} and
  \ref{ex-evl} that the gluing weight times the primitive vector in direction
  of the ray corresponding to $ \alpha $ equals $ H_\alpha \, u_\alpha $, where
  $ H_\alpha $ denotes the Hurwitz number of $ V_1 $ for type I, and the
  product of the Hurwitz numbers of $ V_1 $ and $ V_2 $ for type II. Hence we
  have to show that $ \sum_\alpha H_\alpha \, u_\alpha = 0 $.

  By lemma \ref{lem-forget}, it suffices to prove that $ \sum_\alpha H_\alpha
  \, \ft_I(u_\alpha) = 0 $ for all four-element subsets $ I = \{i,j,k,l\} $ of
  $ I_V $. The combinatorial types $ \alpha $ for which $ \ft_I(u_\alpha) $ is
  a multiple of $ v_{\{i,j\}} $ are exactly the ones corresponding to the three
  cases in lemma \ref{lem-boundary}. Due to the definition of $ u_\alpha $,
  this multiple is $1$, $ d_1 $, and $ d_1+d_2 $, respectively, and hence
  always equal to $ \ord_\calC \ft_I^* (ij|kl) $ for a stable map $ \calC $ of
  this type. As the number of such stable maps is exactly $ H_\alpha $, it
  follows that $ \sum_\alpha H_\alpha \, \ft_I(u_\alpha) $ contains the vector
  $ v_{\{i,j\}} $ with a factor of $ \deg \ft^*_I (ij|kl) $. But the same holds
  for the other two splittings of $I$, and thus we conclude as desired that
    \[ \sum_\alpha H_\alpha \, \ft_I(u_\alpha)
       = \deg \ft^*_I (ij|kl) \, v_{\{i,j\}} +
         \deg \ft^*_I (ik|jl) \, v_{\{i,k\}} +
         \deg \ft^*_I (il|jk) \, v_{\{i,l\}} = 0 \]
  since these three divisors are linearly equivalent and $ v_{\{i,j\}} +
  v_{\{i,k\}} + v_{\{i,l\}} = 0 $ in $ \calM_{0,I} $.
\end{proof}

\begin{proof}[Proof of theorem \ref{thm1}:]
  Theorem \ref{thm-local-balancing} together with the arguments at the
  beginning of section \ref{sec-1dim} shows that all vertices $V$ of
  combinatorial types of $\mathcal{M}_{0,n}(L,\Sigma)$ with $\rdim(V)=1$ are
  good. By lemma \ref{lem-good1dim} we conclude that all vertices are good.
  Hence $\mathcal{M}_{0,n}(L,\Sigma)$ is a tropical variety by theorem
  \ref{thm-glue}, with the weights given in constructions \ref{con-weights} and
  \ref{con:gluing}. The claim about the dimension follows from subsection
  \ref{subsec-ascomplex}.
\end{proof}

\begin{remark}

  By lemma \ref{lem-good1dim}, the case of one-dimensional moduli spaces of
  tropical stable maps to a curve represents a main building block for the
  proof of Theorem \ref{thm1} stating that arbitrary-dimensional moduli spaces
  of tropical stable maps to a curve are balanced. It was also a natural
  starting point for the investigation of the balancing condition for tropical
  moduli spaces of stable maps to a curve. In collaboration with Simon Hampe,
  the polymake extension a-tint \cite{polymake, Ham12} was used to compute ---
  for a large series of relevant examples --- the generating vectors of rays
  for such one-dimensional moduli spaces. GAP \cite{gap4} was used to compute
  conjectural local weights in terms of Hurwitz numbers, and to check the
  balancing condition. These experiments with one-dimensional moduli spaces
  helped us to form a precise conjecture for the weights. Finally, the
  computation of a series of one-dimensional balanced examples led to the proof
  of the balancing condition in the one-dimensional case, and thus also in the
  general case. This work thus heavily relies on the examples computed with the
  help of a-tint and GAP.
\end{remark}

%% file: pics/res1.tex
\begin{picture}(0,0)%
\includegraphics{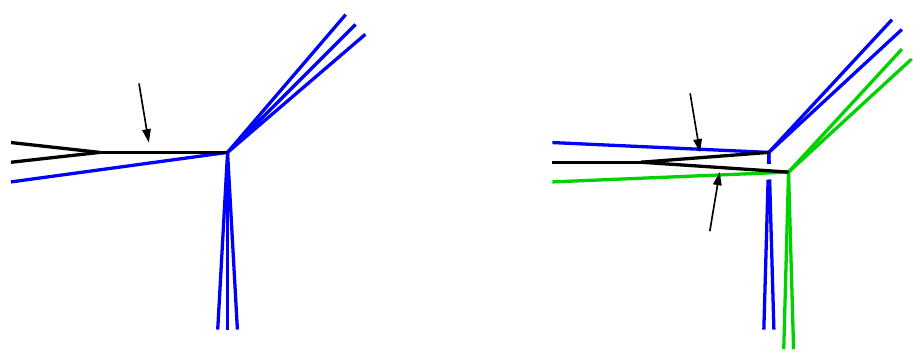}%
\end{picture}%
\setlength{\unitlength}{4144sp}%
\begingroup\makeatletter\ifx\SetFigFont\undefined%
\gdef\SetFigFont#1#2#3#4#5{%
  \reset@font\fontsize{#1}{#2pt}%
  \fontfamily{#3}\fontseries{#4}\fontshape{#5}%
  \selectfont}%
\fi\endgroup%
\begin{picture}(4191,1609)(310,-983)
\put(1756,479){\makebox(0,0)[rb]{\smash{{\SetFigFont{10}{12.0}{\familydefault}{\mddefault}{\updefault}{\color[rgb]{0,0,1}$I_1$}%
}}}}
\put(4276,479){\makebox(0,0)[rb]{\smash{{\SetFigFont{10}{12.0}{\familydefault}{\mddefault}{\updefault}{\color[rgb]{0,0,1}$I_1$}%
}}}}
\put(1351,  6){\makebox(0,0)[rb]{\smash{{\SetFigFont{10}{12.0}{\familydefault}{\mddefault}{\updefault}{\color[rgb]{0,0,1}$V_1$}%
}}}}
\put(3826,  6){\makebox(0,0)[rb]{\smash{{\SetFigFont{10}{12.0}{\familydefault}{\mddefault}{\updefault}{\color[rgb]{0,0,1}$V_1$}%
}}}}
\put(4411,164){\makebox(0,0)[lb]{\smash{{\SetFigFont{10}{12.0}{\familydefault}{\mddefault}{\updefault}{\color[rgb]{0,.82,0}$I_2$}%
}}}}
\put(3961,-286){\makebox(0,0)[lb]{\smash{{\SetFigFont{10}{12.0}{\familydefault}{\mddefault}{\updefault}{\color[rgb]{0,.82,0}$V_2$}%
}}}}
\put(766,-219){\makebox(0,0)[b]{\smash{{\SetFigFont{10}{12.0}{\familydefault}{\mddefault}{\updefault}{\color[rgb]{0,0,0}$V_0$}%
}}}}
\put(3241,-264){\makebox(0,0)[b]{\smash{{\SetFigFont{10}{12.0}{\familydefault}{\mddefault}{\updefault}{\color[rgb]{0,0,0}$V_0$}%
}}}}
\put(3466,276){\makebox(0,0)[b]{\smash{{\SetFigFont{10}{12.0}{\familydefault}{\mddefault}{\updefault}{\color[rgb]{0,0,0}$d_1$}%
}}}}
\put(3556,-579){\makebox(0,0)[b]{\smash{{\SetFigFont{10}{12.0}{\familydefault}{\mddefault}{\updefault}{\color[rgb]{0,0,0}$d_2$}%
}}}}
\put(361,-916){\makebox(0,0)[lb]{\smash{{\SetFigFont{10}{12.0}{\familydefault}{\mddefault}{\updefault}{\color[rgb]{0,0,0}Type I}%
}}}}
\put(2836,-916){\makebox(0,0)[lb]{\smash{{\SetFigFont{10}{12.0}{\familydefault}{\mddefault}{\updefault}{\color[rgb]{0,0,0}Type II}%
}}}}
\put(946,321){\makebox(0,0)[b]{\smash{{\SetFigFont{10}{12.0}{\familydefault}{\mddefault}{\updefault}{\color[rgb]{0,0,0}$d_1$}%
}}}}
\put(325,-138){\makebox(0,0)[rb]{\smash{{\SetFigFont{8}{9.6}{\familydefault}{\mddefault}{\updefault}{\color[rgb]{0,0,0}$x_j$}%
}}}}
\put(325,-25){\makebox(0,0)[rb]{\smash{{\SetFigFont{8}{9.6}{\familydefault}{\mddefault}{\updefault}{\color[rgb]{0,0,0}$x_i$}%
}}}}
\put(2800,-129){\makebox(0,0)[rb]{\smash{{\SetFigFont{8}{9.6}{\familydefault}{\mddefault}{\updefault}{\color[rgb]{0,0,0}$x_i$}%
}}}}
\end{picture}%

%% file: pics/res2.tex
\begin{picture}(0,0)%
\includegraphics{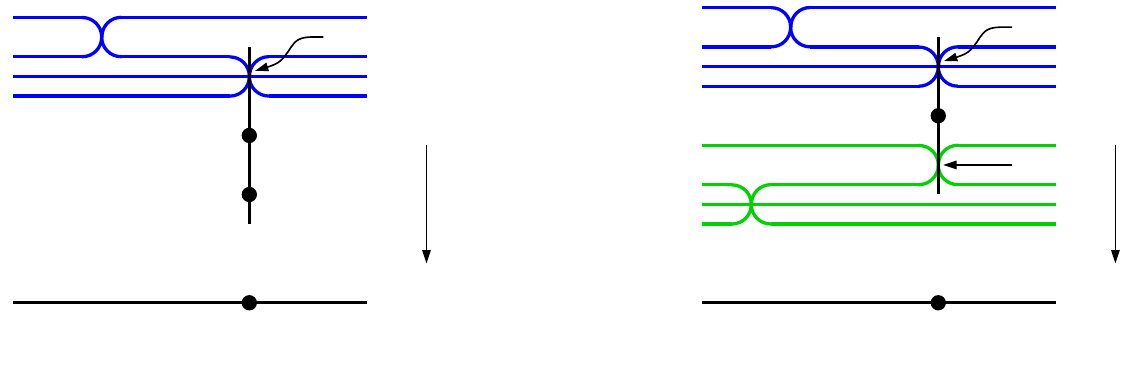}%
\end{picture}%
\setlength{\unitlength}{4144sp}%
\begingroup\makeatletter\ifx\SetFigFont\undefined%
\gdef\SetFigFont#1#2#3#4#5{%
  \reset@font\fontsize{#1}{#2pt}%
  \fontfamily{#3}\fontseries{#4}\fontshape{#5}%
  \selectfont}%
\fi\endgroup%
\begin{picture}(5127,1702)(301,-3001)
\put(361,-2941){\makebox(0,0)[lb]{\smash{{\SetFigFont{10}{12.0}{\familydefault}{\mddefault}{\updefault}{\color[rgb]{0,0,0}Type I}%
}}}}
\put(1509,-1951){\makebox(0,0)[lb]{\smash{{\SetFigFont{10}{12.0}{\familydefault}{\mddefault}{\updefault}{\color[rgb]{0,0,0}$x_i$}%
}}}}
\put(1509,-2221){\makebox(0,0)[lb]{\smash{{\SetFigFont{10}{12.0}{\familydefault}{\mddefault}{\updefault}{\color[rgb]{0,0,0}$x_j$}%
}}}}
\put(1801,-1501){\makebox(0,0)[lb]{\smash{{\SetFigFont{8}{9.6}{\familydefault}{\mddefault}{\updefault}{\color[rgb]{0,0,0}$d_1$}%
}}}}
\put(316,-1591){\makebox(0,0)[rb]{\smash{{\SetFigFont{10}{12.0}{\familydefault}{\mddefault}{\updefault}{\color[rgb]{0,0,1}$C_1$}%
}}}}
\put(3511,-2941){\makebox(0,0)[lb]{\smash{{\SetFigFont{10}{12.0}{\familydefault}{\mddefault}{\updefault}{\color[rgb]{0,0,0}Type II}%
}}}}
\put(4659,-1861){\makebox(0,0)[lb]{\smash{{\SetFigFont{10}{12.0}{\familydefault}{\mddefault}{\updefault}{\color[rgb]{0,0,0}$x_i$}%
}}}}
\put(4951,-1456){\makebox(0,0)[lb]{\smash{{\SetFigFont{8}{9.6}{\familydefault}{\mddefault}{\updefault}{\color[rgb]{0,0,0}$d_1$}%
}}}}
\put(3466,-1546){\makebox(0,0)[rb]{\smash{{\SetFigFont{10}{12.0}{\familydefault}{\mddefault}{\updefault}{\color[rgb]{0,0,1}$C_1$}%
}}}}
\put(3466,-2176){\makebox(0,0)[rb]{\smash{{\SetFigFont{10}{12.0}{\familydefault}{\mddefault}{\updefault}{\color[rgb]{0,.82,0}$C_2$}%
}}}}
\put(2251,-1861){\makebox(0,0)[b]{\smash{{\SetFigFont{10}{12.0}{\familydefault}{\mddefault}{\updefault}{\color[rgb]{0,0,0}$C$}%
}}}}
\put(2251,-2716){\makebox(0,0)[b]{\smash{{\SetFigFont{10}{12.0}{\familydefault}{\mddefault}{\updefault}{\color[rgb]{0,0,0}$\PP^1$}%
}}}}
\put(5401,-1861){\makebox(0,0)[b]{\smash{{\SetFigFont{10}{12.0}{\familydefault}{\mddefault}{\updefault}{\color[rgb]{0,0,0}$C$}%
}}}}
\put(5401,-2716){\makebox(0,0)[b]{\smash{{\SetFigFont{10}{12.0}{\familydefault}{\mddefault}{\updefault}{\color[rgb]{0,0,0}$\PP^1$}%
}}}}
\put(1441,-2851){\makebox(0,0)[b]{\smash{{\SetFigFont{10}{12.0}{\familydefault}{\mddefault}{\updefault}{\color[rgb]{0,0,0}$P_s$}%
}}}}
\put(4591,-2851){\makebox(0,0)[b]{\smash{{\SetFigFont{10}{12.0}{\familydefault}{\mddefault}{\updefault}{\color[rgb]{0,0,0}$P_s$}%
}}}}
\put(4951,-2086){\makebox(0,0)[lb]{\smash{{\SetFigFont{8}{9.6}{\familydefault}{\mddefault}{\updefault}{\color[rgb]{0,0,0}$d_2$}%
}}}}
\put(1419,-2379){\makebox(0,0)[rb]{\smash{{\SetFigFont{10}{12.0}{\familydefault}{\mddefault}{\updefault}{\color[rgb]{0,0,0}$C_0$}%
}}}}
\put(4636,-2244){\makebox(0,0)[lb]{\smash{{\SetFigFont{10}{12.0}{\familydefault}{\mddefault}{\updefault}{\color[rgb]{0,0,0}$C_0$}%
}}}}
\end{picture}%

%% file: pics/res3.tex
\begin{picture}(0,0)%
\includegraphics{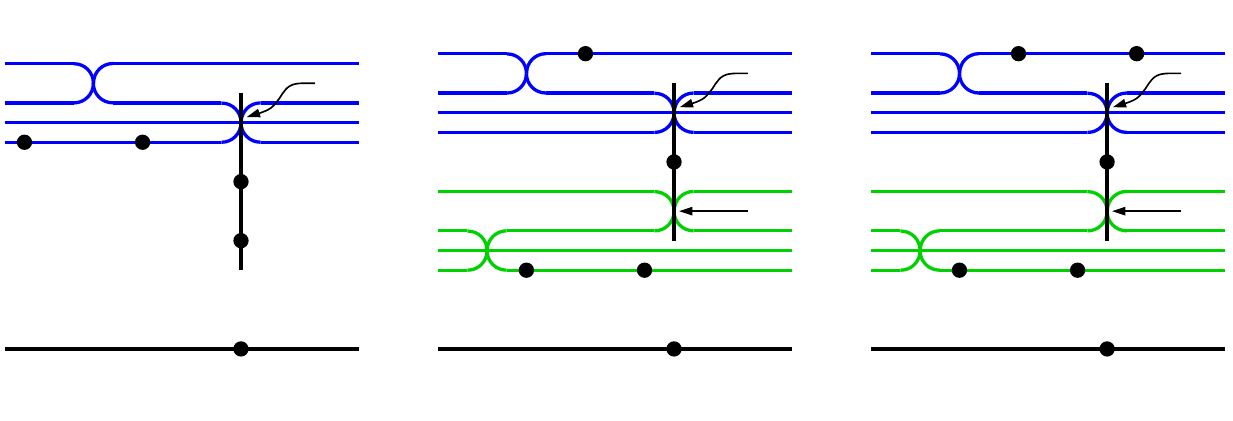}%
\end{picture}%
\setlength{\unitlength}{4144sp}%
\begingroup\makeatletter\ifx\SetFigFont\undefined%
\gdef\SetFigFont#1#2#3#4#5{%
  \reset@font\fontsize{#1}{#2pt}%
  \fontfamily{#3}\fontseries{#4}\fontshape{#5}%
  \selectfont}%
\fi\endgroup%
\begin{picture}(5624,1912)(339,-2996)
\put(1509,-1951){\makebox(0,0)[lb]{\smash{{\SetFigFont{10}{12.0}{\familydefault}{\mddefault}{\updefault}{\color[rgb]{0,0,0}$x_i$}%
}}}}
\put(1509,-2221){\makebox(0,0)[lb]{\smash{{\SetFigFont{10}{12.0}{\familydefault}{\mddefault}{\updefault}{\color[rgb]{0,0,0}$x_j$}%
}}}}
\put(1801,-1501){\makebox(0,0)[lb]{\smash{{\SetFigFont{8}{9.6}{\familydefault}{\mddefault}{\updefault}{\color[rgb]{0,0,0}$d_1$}%
}}}}
\put(1441,-2851){\makebox(0,0)[b]{\smash{{\SetFigFont{10}{12.0}{\familydefault}{\mddefault}{\updefault}{\color[rgb]{0,0,0}$P_s$}%
}}}}
\put(3489,-1861){\makebox(0,0)[lb]{\smash{{\SetFigFont{10}{12.0}{\familydefault}{\mddefault}{\updefault}{\color[rgb]{0,0,0}$x_i$}%
}}}}
\put(3781,-1456){\makebox(0,0)[lb]{\smash{{\SetFigFont{8}{9.6}{\familydefault}{\mddefault}{\updefault}{\color[rgb]{0,0,0}$d_1$}%
}}}}
\put(3421,-2851){\makebox(0,0)[b]{\smash{{\SetFigFont{10}{12.0}{\familydefault}{\mddefault}{\updefault}{\color[rgb]{0,0,0}$P_s$}%
}}}}
\put(5761,-1456){\makebox(0,0)[lb]{\smash{{\SetFigFont{8}{9.6}{\familydefault}{\mddefault}{\updefault}{\color[rgb]{0,0,0}$d_1$}%
}}}}
\put(5401,-2851){\makebox(0,0)[b]{\smash{{\SetFigFont{10}{12.0}{\familydefault}{\mddefault}{\updefault}{\color[rgb]{0,0,0}$P_s$}%
}}}}
\put(361,-2941){\makebox(0,0)[lb]{\smash{{\SetFigFont{10}{12.0}{\familydefault}{\mddefault}{\updefault}{\color[rgb]{0,0,0}(a)}%
}}}}
\put(2341,-2941){\makebox(0,0)[lb]{\smash{{\SetFigFont{10}{12.0}{\familydefault}{\mddefault}{\updefault}{\color[rgb]{0,0,0}(b)}%
}}}}
\put(4321,-2941){\makebox(0,0)[lb]{\smash{{\SetFigFont{10}{12.0}{\familydefault}{\mddefault}{\updefault}{\color[rgb]{0,0,0}(c)}%
}}}}
\put(991,-1861){\makebox(0,0)[b]{\smash{{\SetFigFont{10}{12.0}{\familydefault}{\mddefault}{\updefault}{\color[rgb]{0,0,0}$x_l$}%
}}}}
\put(451,-1861){\makebox(0,0)[b]{\smash{{\SetFigFont{10}{12.0}{\familydefault}{\mddefault}{\updefault}{\color[rgb]{0,0,0}$x_k$}%
}}}}
\put(3286,-2446){\makebox(0,0)[b]{\smash{{\SetFigFont{10}{12.0}{\familydefault}{\mddefault}{\updefault}{\color[rgb]{0,0,0}$x_l$}%
}}}}
\put(2746,-2446){\makebox(0,0)[b]{\smash{{\SetFigFont{10}{12.0}{\familydefault}{\mddefault}{\updefault}{\color[rgb]{0,0,0}$x_k$}%
}}}}
\put(3016,-1231){\makebox(0,0)[b]{\smash{{\SetFigFont{10}{12.0}{\familydefault}{\mddefault}{\updefault}{\color[rgb]{0,0,0}$x_j$}%
}}}}
\put(4726,-2446){\makebox(0,0)[b]{\smash{{\SetFigFont{10}{12.0}{\familydefault}{\mddefault}{\updefault}{\color[rgb]{0,0,0}$x_k$}%
}}}}
\put(5266,-2446){\makebox(0,0)[b]{\smash{{\SetFigFont{10}{12.0}{\familydefault}{\mddefault}{\updefault}{\color[rgb]{0,0,0}$x_l$}%
}}}}
\put(4996,-1231){\makebox(0,0)[b]{\smash{{\SetFigFont{10}{12.0}{\familydefault}{\mddefault}{\updefault}{\color[rgb]{0,0,0}$x_j$}%
}}}}
\put(5536,-1231){\makebox(0,0)[b]{\smash{{\SetFigFont{10}{12.0}{\familydefault}{\mddefault}{\updefault}{\color[rgb]{0,0,0}$x_i$}%
}}}}
\put(3781,-2086){\makebox(0,0)[lb]{\smash{{\SetFigFont{8}{9.6}{\familydefault}{\mddefault}{\updefault}{\color[rgb]{0,0,0}$d_2$}%
}}}}
\put(5761,-2086){\makebox(0,0)[lb]{\smash{{\SetFigFont{8}{9.6}{\familydefault}{\mddefault}{\updefault}{\color[rgb]{0,0,0}$d_2$}%
}}}}
\put(5469,-1861){\makebox(0,0)[lb]{\smash{{\SetFigFont{10}{12.0}{\familydefault}{\mddefault}{\updefault}{\color[rgb]{0,0,0}$x_p$}%
}}}}
\end{picture}%

%% file: pics/res4.tex
\begin{picture}(0,0)%
\includegraphics{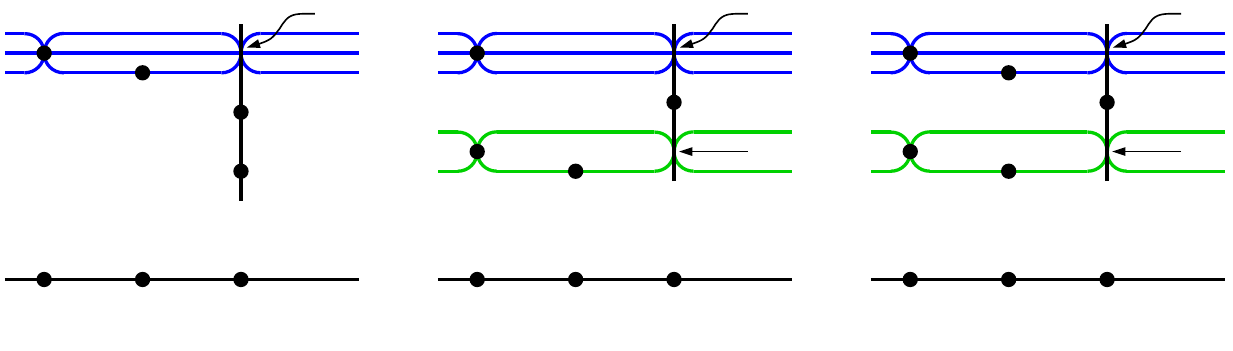}%
\end{picture}%
\setlength{\unitlength}{4144sp}%
\begingroup\makeatletter\ifx\SetFigFont\undefined%
\gdef\SetFigFont#1#2#3#4#5{%
  \reset@font\fontsize{#1}{#2pt}%
  \fontfamily{#3}\fontseries{#4}\fontshape{#5}%
  \selectfont}%
\fi\endgroup%
\begin{picture}(5624,1594)(339,-2996)
\put(1509,-1951){\makebox(0,0)[lb]{\smash{{\SetFigFont{10}{12.0}{\familydefault}{\mddefault}{\updefault}{\color[rgb]{0,0,0}$x_i$}%
}}}}
\put(541,-1771){\makebox(0,0)[b]{\smash{{\SetFigFont{10}{12.0}{\familydefault}{\mddefault}{\updefault}{\color[rgb]{0,0,0}$x_k$}%
}}}}
\put(991,-1861){\makebox(0,0)[b]{\smash{{\SetFigFont{10}{12.0}{\familydefault}{\mddefault}{\updefault}{\color[rgb]{0,0,0}$x_l$}%
}}}}
\put(1509,-2221){\makebox(0,0)[lb]{\smash{{\SetFigFont{10}{12.0}{\familydefault}{\mddefault}{\updefault}{\color[rgb]{0,0,0}$x_j$}%
}}}}
\put(1801,-1501){\makebox(0,0)[lb]{\smash{{\SetFigFont{8}{9.6}{\familydefault}{\mddefault}{\updefault}{\color[rgb]{0,0,0}$d_1$}%
}}}}
\put(1441,-2851){\makebox(0,0)[b]{\smash{{\SetFigFont{10}{12.0}{\familydefault}{\mddefault}{\updefault}{\color[rgb]{0,0,0}$0$}%
}}}}
\put(991,-2851){\makebox(0,0)[b]{\smash{{\SetFigFont{10}{12.0}{\familydefault}{\mddefault}{\updefault}{\color[rgb]{0,0,0}$1$}%
}}}}
\put(541,-2851){\makebox(0,0)[b]{\smash{{\SetFigFont{10}{12.0}{\familydefault}{\mddefault}{\updefault}{\color[rgb]{0,0,0}$\infty$}%
}}}}
\put(1981,-2941){\makebox(0,0)[rb]{\smash{{\SetFigFont{10}{12.0}{\familydefault}{\mddefault}{\updefault}{\color[rgb]{0,0,0}(a)}%
}}}}
\put(3961,-2941){\makebox(0,0)[rb]{\smash{{\SetFigFont{10}{12.0}{\familydefault}{\mddefault}{\updefault}{\color[rgb]{0,0,0}(b)}%
}}}}
\put(5941,-2941){\makebox(0,0)[rb]{\smash{{\SetFigFont{10}{12.0}{\familydefault}{\mddefault}{\updefault}{\color[rgb]{0,0,0}(c)}%
}}}}
\put(3421,-2851){\makebox(0,0)[b]{\smash{{\SetFigFont{10}{12.0}{\familydefault}{\mddefault}{\updefault}{\color[rgb]{0,0,0}$\infty$}%
}}}}
\put(2521,-2851){\makebox(0,0)[b]{\smash{{\SetFigFont{10}{12.0}{\familydefault}{\mddefault}{\updefault}{\color[rgb]{0,0,0}$0$}%
}}}}
\put(2971,-2851){\makebox(0,0)[b]{\smash{{\SetFigFont{10}{12.0}{\familydefault}{\mddefault}{\updefault}{\color[rgb]{0,0,0}$1$}%
}}}}
\put(3781,-1501){\makebox(0,0)[lb]{\smash{{\SetFigFont{8}{9.6}{\familydefault}{\mddefault}{\updefault}{\color[rgb]{0,0,0}$d_1$}%
}}}}
\put(3781,-2131){\makebox(0,0)[lb]{\smash{{\SetFigFont{8}{9.6}{\familydefault}{\mddefault}{\updefault}{\color[rgb]{0,0,0}$d_2$}%
}}}}
\put(3489,-1906){\makebox(0,0)[lb]{\smash{{\SetFigFont{10}{12.0}{\familydefault}{\mddefault}{\updefault}{\color[rgb]{0,0,0}$x_i$}%
}}}}
\put(2521,-1771){\makebox(0,0)[b]{\smash{{\SetFigFont{10}{12.0}{\familydefault}{\mddefault}{\updefault}{\color[rgb]{0,0,0}$x_j$}%
}}}}
\put(2521,-2221){\makebox(0,0)[b]{\smash{{\SetFigFont{10}{12.0}{\familydefault}{\mddefault}{\updefault}{\color[rgb]{0,0,0}$x_k$}%
}}}}
\put(2971,-2311){\makebox(0,0)[b]{\smash{{\SetFigFont{10}{12.0}{\familydefault}{\mddefault}{\updefault}{\color[rgb]{0,0,0}$x_l$}%
}}}}
\put(5761,-1501){\makebox(0,0)[lb]{\smash{{\SetFigFont{8}{9.6}{\familydefault}{\mddefault}{\updefault}{\color[rgb]{0,0,0}$d_1$}%
}}}}
\put(5761,-2131){\makebox(0,0)[lb]{\smash{{\SetFigFont{8}{9.6}{\familydefault}{\mddefault}{\updefault}{\color[rgb]{0,0,0}$d_2$}%
}}}}
\put(5401,-2851){\makebox(0,0)[b]{\smash{{\SetFigFont{10}{12.0}{\familydefault}{\mddefault}{\updefault}{\color[rgb]{0,0,0}$\infty$}%
}}}}
\put(4951,-2851){\makebox(0,0)[b]{\smash{{\SetFigFont{10}{12.0}{\familydefault}{\mddefault}{\updefault}{\color[rgb]{0,0,0}$1$}%
}}}}
\put(4501,-2851){\makebox(0,0)[b]{\smash{{\SetFigFont{10}{12.0}{\familydefault}{\mddefault}{\updefault}{\color[rgb]{0,0,0}$0$}%
}}}}
\put(4501,-1771){\makebox(0,0)[b]{\smash{{\SetFigFont{10}{12.0}{\familydefault}{\mddefault}{\updefault}{\color[rgb]{0,0,0}$x_i$}%
}}}}
\put(4501,-2221){\makebox(0,0)[b]{\smash{{\SetFigFont{10}{12.0}{\familydefault}{\mddefault}{\updefault}{\color[rgb]{0,0,0}$x_l$}%
}}}}
\put(4951,-1861){\makebox(0,0)[b]{\smash{{\SetFigFont{10}{12.0}{\familydefault}{\mddefault}{\updefault}{\color[rgb]{0,0,0}$x_j$}%
}}}}
\put(4951,-2311){\makebox(0,0)[b]{\smash{{\SetFigFont{10}{12.0}{\familydefault}{\mddefault}{\updefault}{\color[rgb]{0,0,0}$x_k$}%
}}}}
\put(5469,-1906){\makebox(0,0)[lb]{\smash{{\SetFigFont{10}{12.0}{\familydefault}{\mddefault}{\updefault}{\color[rgb]{0,0,0}$x_p$}%
}}}}
\end{picture}%

%% file: main.bbl
\providecommand{\bysame}{\leavevmode\hbox to3em{\hrulefill}\thinspace}
\providecommand{\MR}{\relax\ifhmode\unskip\space\fi MR }
\providecommand{\MRhref}[2]{%
  \href{http://www.ams.org/mathscinet-getitem?mr=#1}{#2}
}
\providecommand{\href}[2]{#2}
\begin{thebibliography}{GKM09}

\bibitem[All12]{All09}
Lars Allermann, \emph{Tropical intersection products on smooth varietie},
  Journal of the EMS \textbf{14} (2012), no.~1, 107--126.

\bibitem[AR10]{AR07}
Lars Allermann and Johannes Rau, \emph{First steps in tropical intersection
  theory}, Math. Z. \textbf{264} (2010), no.~3, 633--670, arXiv:0709.3705.

\bibitem[BBM11]{BBM10}
Beno\^{\i}t Bertrand, Erwan Brugall\'{e}, and Grigory Mikhalkin,
  \emph{{Tropical Open Hurwitz numbers}}, Rend. Semin. Mat. Univ. Padova
  \textbf{125} (2011), 157--171.

\bibitem[BBM14]{BBM11}
\bysame, \emph{Genus 0 characteristic numbers of tropical projective plane},
  Compos. Math. \textbf{150} (2014), no.~1, 46--104, arXiv:1105.2004.

\bibitem[Beh97]{Behrend}
K.~Behrend, \emph{Gromov-{W}itten invariants in algebraic geometry}, Invent.
  Math. \textbf{127} (1997), no.~3, 601--617.

\bibitem[BF97]{BF97}
Kai Behrend and Barbara Fantechi, \emph{The intrinsic normal cone}, Invent.
  Math. \textbf{128} (1997), 45--88.

\bibitem[BM13a]{BruMa}
Erwan Brugall\'{e} and Hannah Markwig, \emph{Deformation of tropical
  {Hirzebruch} surfaces and enumerative geometry}, {Preprint}, arXiv:1303.1340,
  J. Algebraic Geom. (to appear), 2013.

\bibitem[BM13b]{BM13}
Arne Buchholz and Hannah Markwig, \emph{Tropical covers of curves and their
  moduli spaces}, Comm. Contemp. Math. (2013), doi:10.1142/S0219199713500454.

\bibitem[Cap14]{Cap12}
Lucia Caporaso, \emph{Gonality of algebraic curves and graphs}, Springer Proc.
  Math. Stat. \textbf{71} (2014), 77--108.

\bibitem[CJM10]{CJM10}
Renzo Cavalieri, Paul Johnson, and Hannah Markwig, \emph{{Tropical Hurwitz
  numbers}}, J. Algebr. Comb. \textbf{32} (2010), no.~2, 241--265,
  arXiv:0804.0579.

\bibitem[CJM11]{cjm:wcfdhn}
Renzo Cavalieri, Paul Johnson, and Hannah Markwig, \emph{Wall crossings for
  double {H}urwitz numbers}, Adv. Math. \textbf{228} (2011), no.~4, 1894--1937,
  arXiv:1003.1805.

\bibitem[CMR14]{CMR14}
Renzo Cavalieri, Hannah Markwig, and Dhruv Ranganathan, \emph{Tropicalizing the
  space of admissible covers}, {Preprint}, arXiv:1401.4626, Math. Ann. (to
  appear), 2014.

\bibitem[CMR16]{CMR14b}
\bysame, \emph{Tropical compactification and the {Gromov--Witten} theory of
  $\mathbb{P}^1$}, Selecta Math. \textbf{1-34} (2016),
  doi:10.1007/s00029-016-0265-7, arXiv:1410.2837.

\bibitem[FP97]{FP97}
William Fulton and Rahul Pandharipande, \emph{Notes on stable maps and quantum
  cohomology}, Algebraic Geometry, Santa Cruz 1995 (J{\'a}nos~Koll{\'a}r
  et~al., ed.), Proceedings of Symposia in Pure Mathematics, no. 62,2, Amer.\
  Math.\ Soc., 1997, pp.~45--96.

\bibitem[FS05]{FS05}
Eva~Maria Feichtner and Bernd Sturmfels, \emph{Matroid polytopes, nested sets
  and {Bergman} fans}, Portugaliae Mathematica \textbf{62} (2005), 437--468,
  arXiv:math.CO:0411260.

\bibitem[GAP16]{gap4}
The GAP~Group, \emph{{GAP -- Groups, Algorithms, and Programming, Version
  4.8.6}}, 2016.

\bibitem[GJ00]{polymake}
Ewgenij Gawrilow and Michael Joswig, \emph{polymake: a framework for analyzing
  convex polytopes}, Polytopes --- Combinatorics and Computation (Gil Kalai and
  G\"unter~M. Ziegler, eds.), Birkh\"auser, 2000, pp.~43--74.

\bibitem[GKM09]{GKM07}
Andreas Gathmann, Michael Kerber, and Hannah Markwig, \emph{Tropical fans and
  the moduli space of rational tropical curves}, Compos. Math. \textbf{145}
  (2009), no.~1, 173--195, arXiv:0708.2268.

\bibitem[GO17]{GO14}
Andreas Gathmann and Dennis Ochse, \emph{Moduli spaces of tropical curves in
  tropical varieties}, arXiv:1705.07626, 2017.

\bibitem[Gro14]{Gro14}
Andreas Gross, \emph{Correspondence theorems via tropicalizations of moduli
  spaces}, Comm. Contemp. Math. (to appear), arXiv:1406.1999, 2014.

\bibitem[Ham14]{Ham12}
Simon Hampe, \emph{a-tint: a polymake extension for algorithmic tropical
  intersection theory}, European J. Combin. \textbf{36C} (2014), 579--607,
  arXiv:1208.4248.

\bibitem[Mik05]{Mi03}
Grigory Mikhalkin, \emph{Enumerative tropical geometry in {${\mathbb{R}^2}$}},
  J. Amer. Math. Soc. \textbf{18} (2005), 313--377, arXiv:math.AG/0312530.

\bibitem[Mik07]{Mi07}
\bysame, \emph{Moduli spaces of rational tropical curves}, Proceedings of
  G\"okova Geometry-Topology Conference GGT 2006 (2007), 39--51,
  arXiv:0704.0839.

\bibitem[MR09]{MR08}
Hannah Markwig and Johannes Rau, \emph{{Tropical descendant Gromov-Witten
  invariants}}, Manuscripta Math. \textbf{129} (2009), no.~3, 293--335,
  arXiv:0809.1102.

\bibitem[Och13]{O13}
Dennis Ochse, \emph{Moduli spaces of rational tropical stable maps into smooth
  tropical varieties}, Ph.D. thesis, TU Kaiserslautern, 2013.

\bibitem[SS04]{SS04a}
David Speyer and Bernd Sturmfels, \emph{The tropical {Grassmannian}}, Adv.
  Geom. \textbf{4} (2004), 389--411.

\bibitem[Vak00]{Vak00}
Ravi Vakil, \emph{The enumerative geometry of rational and elliptic curves in
  projective space}, J. Reine Angew. Math. (Crelle's Journal) \textbf{529}
  (2000), 101--153.

\end{thebibliography}
